\documentclass[12pt]{article}
\usepackage[flushleft]{threeparttable}
\usepackage{amssymb}
\usepackage{amsmath}
\usepackage{graphicx}
\usepackage{tabularx}
\usepackage{calc}
\usepackage{tikz}
\usepackage{float}
\usepackage{amsmath}
\usepackage{amsthm}
\usepackage{longtable}
\usetikzlibrary{decorations.markings}
\tikzstyle{vertex}=[circle, draw, inner sep=0pt, minimum size=6pt]

\newtheorem{theorem}{Theorem}[section]

\newtheorem{lemma}[theorem]{Lemma}

\newtheorem{conjecture}[theorem]{Conjecture}

\newtheorem{construction}[theorem]{Construction}
\newtheorem{proposition}[theorem]{Proposition}

\newenvironment{pf}           {\noindent{\bf Proof.} }%
                                {\null\hfill$\Box$\par\medskip\medskip\medskip\medskip}

\begin{document}

%%%%%%%%%%%%%%%%%%%%%%%%%%%
%%  Article Begins       %%
%%%%%%%%%%%%%%%%%%%%%%%%%%%
\title{Graceful Labellings of Variable Windmills Using Skolem Sequences}

\author{Ahmad H. Alkasasbeh \hspace{0.5in} Danny Dyer\\
Department of Mathematics and Statistics\\
St. John's Campus, Memorial University of Newfoundland\\
St. John's, Newfoundland\\
Canada\\
{\tt ahmad84@mun.ca}\hspace{0.5in}{\tt dyer@mun.ca}
\and
Jared Howell\\
School of Science and the Environment\\
Grenfell Campus, Memorial University of Newfoundland\\
Corner Brook, Newfoundland\\
Canada\\
{\tt jahowell@grenfell.mun.ca}}
\maketitle
 %\vspace{0.5in}

%Jared Howell \footnote{jahowell@grenfell.mun.ca}\\
%Division of Science (Mathematics)\\
%Grenfell Campus, Memorial University of Newfoundland\\
%Corner Brook, Newfoundland\\
%Canada
\begin{abstract}
In this paper, we introduce graceful and near graceful labellings of several families of windmills. In particular, we use Skolem-type sequences to prove (near) graceful labellings exist for windmills with $C_{3}$ and $C_{4}$ vanes, and infinite families of $3,5$-windmills and $3,6$-windmills. Furthermore, we offer a new solution showing that the graph obtained from the union of $t$ $5$-cycles with one vertex in common ($C_{5}^{t}$) is graceful if and only if $t\equiv0,3\ ({\rm mod}\ 4)$ and near graceful when $t\equiv1,2\ ({\rm mod}\ 4)$.
%In this paper we introduce graceful and near graceful labellings of several families of windmills. In particular, we use Skolem-like sequences to prove (near) graceful labellings exist for infinite families of $3$ and $4$ vane windmills, $3,5$-windmills and $3,6$-windmills. Furthermore, we offer a new solution showing that $G=C_{5}^{t}$ is near graceful if and only if $t\equiv1,2\ ({\rm mod}\ 4)$.
%for the $C_{3}^{(t)}C_{4}^{(s)}$, $C_{3}^{(t)}C_{5}^{(s)}$ and $C_{3}^{(t)}C_{6}^{(s)}$. Furthermore, we verified that $G=C_{5}^{(t)}$ is near graceful if and only if $t\equiv1,2\ ({\rm mod}\ 4)$.
\end{abstract}


%\clearpage
%%%%%%%%%%%%%%%%%%%%%%%%%
%   Intro               %
%%%%%%%%%%%%%%%%%%%%%%%%%
%\tableofcontents
%%%%%%%%%%%%%%%%%%%%%%%%%%%%%%%%%%%%%%%%%%%%%%%%%%%%%%%  Intoduction  %%%%%%%%%%%%%%%%%%%%%%%%%%%%%%%%%%%%%%%%%%%%%%%%%%%%%%%%%%%%%%%%%%%%%%%%%%%%%%%%%%%%%%%%%%%%%%%%%%
%\section{Introduction}\label{section1}
%\chapter{Graceful Labellings of Variable Windmills Using Skolem Sequences}\label{ch:3}
%{\let\thefootnote\relax\footnote{{Submitted to the JCMCC-Journal of Combinatorial Mathematics and Combinatorial Computing.}}}
%%%%%%%%%%%%%%%%%%%%%%%%%%%%%%%%%%%%%%%%%%%%%%%%%%%Introduction%%%%%%%%%%%%%%%%%%%%%%%%%%%%%%%%%%%%%%%%%%%%%%%%%%%
\section{Introduction}\label{section1}
In \cite{rosa1}, Rosa introduced a new type of graph labelling known as a $\beta$-labelling, or graceful labelling, as it was renamed later. Let $G=(V,E)$ be a graph with $m$ edges. Let $f:V(G)\rightarrow \left\{0,1,2,\ldots,m\right\}$ be a labelling of $V$ of $G$ and let $g:E(G)\rightarrow \left\{1,2,\ldots,m\right\}$ be the induced edge labelling defined by $g(uv)=|f(u)-f(v)|,$ for all $uv \in E.$ The labelling $f$ is said to be a \textit{graceful labelling} if and only if $f$ is an injective mapping and $g$ is a bijection. If a graph $G$ has a graceful labelling then we say $G$ is graceful.

A near graceful labelling of a graph $G=(V,E)$ with $m$ edges is defined in a similar way. Let $f:V(G)\rightarrow \left\{0,1,2,\ldots,m+1\right\}$ be a labelling of $V$ of $G$ and let $g:E(G)\rightarrow A$ be the induced edge labelling defined by $g(uv)=|f(u)-f(v)|,$ for all $uv \in E$, where $A$ is $\{1,2,\dots,m-1,m\}$ or $\{1,2,\dots,m-1,m+1\}$. The labelling $f$ is said to be a \textit{near graceful labelling} if and only if $f$ is an injective mapping and $g$ is a bijection. If a graph $G$ has a near graceful labelling then we say $G$ is near graceful. In this paper, all near graceful labellings constructed will omit the vertex label $m$ and the edge label $m$.

In this paper, we adopt the convention that $0$ is a natural number. So, when we write $\left[a,b\right]$ with $a,b \in \mathbb{N}$ and $a<b$, we are indicating the set $\left\{x \in \mathbb{N}| a \leq x \leq b\right\}$.

Bermond, in \cite{Bermond}, proved that Dutch windmills (the graphs consisting of $t$ copies of $K_3$ with one vertex in common) are graceful. Let $C_{n}$ be a cycle of length $n \geq 3$, and $C_{n}^{t}$ be the graph obtained from the union of $t$ $n$-cycles with one vertex in common that we will call the \textit{central vertex}.  In \cite{koh}, the authors stated the following conjecture.

\begin{conjecture}\cite{koh} $C_{n}^{t}$ is graceful if and only if $nt\equiv 0,3\ ($mod$\ 4)$.\end{conjecture}

This conjecture has been shown to hold for $n=3, n=4, n=5,$ and $n=6$ and $n=4k$ ($k$ any positive integer) in \cite{Bermond1}, \cite{kejie}, \cite{yang}, and \cite{koh}, respectively. Also, in \cite{yang7,yang9,yang11,yang13}, the authors show that the graceful labellings exist for the $C_{n}^{t}$ with $n=7,9,11,13$, respectively. In 2012, Dyer, et al. \cite{dyer}, use Skolem sequences to prove that all Dutch windmills with zero, one or two pendant triangles are (near) graceful. A comprehensive survey of graceful labelling can be found in \cite{Gallian}.

%%%%%%%%%%%%%%%%%%%%%%%%%%%%%%%%%%%%%%%%%%%%%%%%%%%%%%%%%%%%%%%%%%%%%%%%%%%%%%%%%%%%%%%%%%%%%%%%%%%%%%%%%%%%%%%%%%%%%%%%%%%%%%%%%%%%%%%%%%%%%%%%%%%%%%%
We define an $m,n$-windmill to be a graph $G=C_{m}^{s}C_{n}^{t}$ obtained from identifying the central vertices of $C_{m}^{s}$ and $C_{n}^{t}$, where $m\neq n$. In other words, the graph $C_{m}^{s}C_{n}^{t}$ is a windmill with $s$ $m$-cycle vanes and $t$ $n$-cycle vanes. More generally, we call any windmill made up of two or more cycle lengths a {\it variable windmill}.

\begin{figure}[H]
\begin{center}
\linethickness{0.7pt}
 \resizebox{0.65\textwidth}{!}{%
\begin{tikzpicture}
%First place the vertices. Note the names correspond to the labels.
\node[circle,draw, minimum size=8mm] (v0) at (0,0) {$\textbf{0}$};
\node[circle,draw, minimum size=8mm] (v18) at (0:5) {$\textbf{18}$};
\node[circle,draw, minimum size=8mm] (v20) at (42:5) {$\textbf{20}$};
\node[circle,draw, minimum size=8mm] (v15) at (105:5) {$\textbf{15}$} ;
\node[circle,draw, minimum size=8mm] (v1) at (129.5:8) {$\textbf{1}$} ;
\node[circle,draw, minimum size=8mm] (v7) at (144:5) {$\textbf{7}$} ;	
\node[circle,draw, minimum size=8mm] (v23) at (54:5) {$\textbf{23}$} ;
\node[circle,draw, minimum size=8mm] (v24) at (93:5) {$\textbf{24}$} ;
\node[circle,draw, minimum size=8mm] (v11) at (156:5) {$\textbf{11}$} ;
\node[circle,draw, minimum size=8mm] (v2) at (175.5:8) {$\textbf{2}$} ;
\node[circle,draw, minimum size=8mm] (v12) at (195:5) {$\textbf{12}$} ;
\node[circle,draw, minimum size=8mm] (v17) at (258:5) {$\textbf{17}$} ;
\node[circle,draw, minimum size=8mm] (v21) at (297:5) {$\textbf{21}$} ;
\node[circle,draw, minimum size=8mm] (v16) at  (207:5) {$\textbf{16}$} ;
\node[circle,draw, minimum size=8mm] (v3) at (226.5:8) {$\textbf{3}$} ;
\node[circle,draw, minimum size=8mm] (v8) at (246:5) {$\textbf{8}$} ;
\node[circle,draw, minimum size=8mm] (v22) at (309:5) {$\textbf{22}$} ;
\node[circle,draw, minimum size=8mm] (v19) at (348:5) {$\textbf{19}$} ;
%Next do the edges and labels. Note the label positions; you can make them more precise by using angles instead of 8 directions right, above, etc.
\foreach \x/\y/\lab in {0/18/above,0/20/below right,0/15/left,0/7/above right,0/23/above left,0/24/right,0/11/below,0/12/above,0/17/right,0/21/left,0/16/below,0/8/above left,0/22/above right,0/19/below,15/1/below,7/1/below right,11/2/below,12/2/above,16/3/right,8/3/above,18/20/left,23/24/below,17/21/above,22/19/above left}
\draw[thick] (v\x)--node[above, \lab, text=blue] {\pgfmathparse{int(abs(\x-\y))}{\bf{\pgfmathresult}}} (v\y);
%\draw[thick] (v\x)--node[above, \lab, text=blue] {\pgfmathparse{int(abs(\x-\y))}{\pgfmathresult}} (v\y);
\end{tikzpicture}
}
\caption{Graceful labelling of $C_{3}^{4}C_{4}^{3}$.}\label{ga} % same our construction section c34
\end{center}
\end{figure}
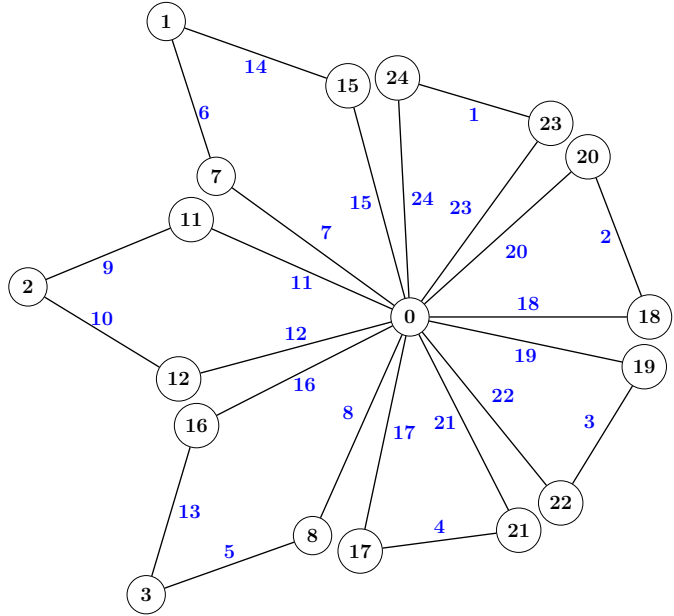

In this paper we use Skolem-type sequences to show (near) graceful labellings exist for the $G=C_{n}^{t}C_{m}^{s}$ graphs where $n=3$ and $m=4,5,6$. An example of a graceful labelling of $C_{3}^{4}C_{4}^{3}$ is given in Figure~\ref{ga}.
%%%%%%%%%%%%%%%%%%%%%%%%%%%%%%%%%%%%%%%%%%%%%%%%%%%%%%%%%%%%%%%%%%%%%%%%%%%%%%%%%%%%%%%%%%%%%%%%%%%%%%%%%%%%%%%%%%%%%%%%%%%%%%%%%%%%%%%%%%%%%%%%%%%%%%%%%%%%%%%%%%%%%%%%
%%%%%%%%%%%%%%%%%%%%%%%%%%%%%%%%%%%%%%%%%%%%%%%%%%%%%%%%%%%%%%%%%%%%%%%%%%%%%%%%%%%%%%%%%%%%%%%%%%%%%%%%%%%%%%%%%%%%%%%%%%%%%%%%%%%%%%%%%%%%%%%%%%%%%%%%%%%%%%%%%%%%%%%%
%%%%%%%%%%%%%%%%%%%%%%%%%%%%%%%%%%%%%%%%%%%%%%%%%%%%%%%%%%%%%%%%%%%%%%%%%%%%%%%%%%%%%%%%%%%%%%%%%%%%%%%%%%%%%%%%%%%%%%%%%%%%%%%%%%%%%%%%%%%%%%%%%%%%%%%%%%%%%%%%%%%%%%%%%%%%%%%%%%%%%%%%%%%%%%%%%%%%%%%%%%%%%%%%%%%%%%%%%%%%%%%%%%%%%%%%%%%%%%%%%%%%%%%%%%%%%%%%%%%%%%%%%%%%%%%%%%%%%%%%%%%%%%%%%%%%%%%%%%%%%%%%%%%%%%%%%%%%%%%%%%%%%%%%%%%%%%%%
%%%%%%%%%%%%%%%%%%%%%%%%%%%%%%%%%%%%%%%%%%%%%%%%%%%%%%%  Preliminaries  %%%%%%%%%%%%%%%%%%%%%%%%%%%%%%%%%%%%%%%%%%%%%%%%%%%%%%%%%%%%%%%%%%%%%%%%%%%%%%%%%%%%%%%%%%%%%%%%
\section{Skolem-type Sequences}\label{section2}
%In this section we will present a new definition which is a Skolem-like sequence and several definitions and known results for Skolem and Langford sequences.\\
In this section we begin by defining a Skolem-type sequence, and then present several definitions and known results for Skolem and Langford sequences.

A {\it Skolem-type sequence} of order $n$ is a sequence $\left(s_{1},s_{2},\ldots,s_{2n}\right)$ of $2n$ integers satisfying the conditions:
\begin{enumerate}
%\item $\forall\ k \in H$, where $H$ is a set of $n$ distinct positive integers, $\exists\ s_{i},s_{j}\in K$ such that $s_{i}=s_{j}=k;$
\item for some set $H$ of $n$ distinct positive integers, $\forall\ k \in H$, $\exists\ s_{i},s_{j}$ such that $s_{i}=s_{j}=k;$
\item if $s_{i}=s_{j}=k,$ with $i < j,$ then $j-i=k$.
\end{enumerate}

For example, $(6,4,1,1,3,4,6,3)$ is a Skolem-type sequence of order $4$ with $H=\{1,3,4,6\}$. We can also write Skolem-type sequences by specifying the ordered pairs where identical elements of $H$ occur, as follows $\{(a_i,b_i): i \in H, b_i-a_i=i\}$. So, we could equivalently write the sequence $\{(3,4),(5,8),(2,6),(1,7)\}$. In the Skolem-type sequences if $b_i > a_i$ we call $b_i$ \textit{right endpoint}. The definition of the Skolem-type sequence introduced in \cite{baker2}.

A {\it hooked Skolem-type sequence} of order $n$ is a sequence $\left(s_{1},s_{2},\ldots,s_{2n+1}\right)$ of $2n+1$ integers satisfying the conditions of a Skolem-type sequence with the added condition that $s_{2n}=0$. For example, $(5,3,1,1,3,5,2,0,2)$ is a hooked Skolem-type sequence of order $4$ with $H=\left\{1,2,3,5\right\}$.
%%%%%%%%%%%%%%%%%%%%%%%%%%%%%%%%%%%%%%%%%%%%%%%%%%%%%%%%%%%%%%%%%%%%%%%%%%%
%\rm A {\it hooked Skolem-like sequence} of order $n$ is a sequence $hK=\left(s_{1},s_{2},\ldots,s_{2n+1}\right)$ of $2n$ integers satisfying these conditions:
%\begin{enumerate}
%\item $\forall\ k \in H$, where $H$ is a set of $n$ distinct positive integers, $\left|H\right|=n$, $\exists\ s_{i},s_{j}\in K$ such that $s_{i}=s_{j}=k;$
%\item if $s_{i}=s_{j}=k,$ with $i < j,$ then $j-i=k$.
%\item $s_{2n}=0.$\left[ \right]
%\end{enumerate}
%For the following different types of sequences, we use the definitions from the Handbook of Combinatorial Designs \cite{crc}.\par
%%%%%%%%%%%%%%%%%%%%%%%%%%%%%%%%%%%%%%%%%%%%%%%%%%%%% Lang %%%%%%%%%%%%%%%%%%%%%%%%%%%%%%%%%%%%%%%%%%%%%%%%%%%%%%%%%%%%%%%%%%%%%%%%%%%%%%%%%%%%%%%%%%%%%%%%%%%%%%%%%%%

A \textit{ (hooked) Skolem sequence} of order $n$ is a (hooked) Skolem-type sequence of order $n$ with $H = \left[1,n\right]$. Necessary and sufficient conditions for the existence of Skolem sequences are given in \cite{skolem}.

% A \textit{hooked Skolem sequence} of order $n$ is a hooked Skolem-like sequence of order $n$ with $H = \left\{1,2,\ldots,n\right\}.$
%%%%%%%%%%%%%%%%%%%%%%%%%%%%%%%%%%%%%%%%%%%%%%%%%%%%%%%%%%%%%%%%%%%%%%%%%%%%%%%%%%%%%%%%%%%%%%%%%%%%%%%%%%%%%%%%%%%%%%%%%%%%%%%%%%%%%%%%%%%%%%%%%%%%%%%%%%%%%%%%%%%%%%%
%%%%%%%%%%%%%%%%%%%%%%%%%%%%%%%%%%%%%%%%%%%%%%%%%%%%% NEAR %%%%%%%%%%%%%%%%%%%%%%%%%%%%%%%%%%%%%%%%%%%%%%%%%%%%%%%%%%%%%%%%%%%%%%%%%%%%%%%%%%%%%%%%%%%%%%%%%%%%%%%%%%%
%In this section, we need to discuss another type of sequences which is a near Skolem sequences and related constructions.\\
Let $m$ and $n$ be positive integers, with $m\leq n$. A {\em (hooked) near-Skolem sequence} of order $n$ and defect $m$ is a sequence $(s_{1},s_{2},\ldots,s_{2n-2})$ is a (hooked) Skolem-type sequence of order $n-1$ with $H=\left[1,m-1\right] \cup \left[m+1,n\right]$. For example, $(1,1,6,3,7,5,3,2,6,2,5,7)$ is a $4$-near-Skolem sequence of order $7$ and $(2,5,2,4,6,7,5,4,1,1,6,0,7)$ is a hooked $3$-near-Skolem sequence of order $7$. Necessary and sufficient conditions for the existence of near-Skolem sequences are given in \cite{shalaby}.

A \textit{(hooked) Langford sequence} with defect $d$ and order $l$ is a (hooked) Skolem-type sequence with $H = \left[d,d+l-1\right]$. Necessary and sufficient conditions for the existence of Langford sequences are given in \cite{simpson}.

Note that in this paper, in the constructions, $a_i$ and $b_i$ or $c_j$ and $d_j$ or $e_j$ and $f_j$ represent the two positions in the Skolem-type sequence of the element $i$ and $j$, with $a_i < b_i$, $c_j < d_j$ and $e_j < f_j$ with $1 \leq i,j \leq n$. When written in this form, we say that $a_i$, $c_j$, and $e_j$ are left endpoints, and $b_i$, $d_j$, and $f_j$ are right endpoints.
%A \textit{hooked Langford sequence} of defect $d$ and order $l$ is a hooked Skolem-like sequence with $H = \left\{d,d+1,\ldots,d+l-1\right\}$ and $\left|H\right|=l$.
%%%%%%%%%%%%%%%%%%%%%%%%%%%%%%%%%%%%%%%%%%%%%%%%%%%%%%%% langford construction %%%%%%%%%%%%%%%%%
\begin{construction}\label{ldd}
From Table~\ref{tablel}, we can construct a Langford sequence with defect $d \geq 1$ and order $2d-1$, (omitting row $2$ when $d=1$). We define $L_{d}^{2d-1}$ to be exactly this sequence for $d\ge 1$.
%\begin{table}[H]
%\begin{center}
%\scalebox{0.85}{
%\begin{tabular}{|r|l|l|}
\begin{table}[ht]
\begin{center}%\setlength\extrarowheight{7pt}\renewcommand{\arraystretch}{0.9}
\scalebox{1.00}{
\begin{tabular}{|c|c|c|c|c|}
\hline
  $i$ & $a_i$ & $b_i$ &  \\
\hline
\hline
 $d+2r$ &\ $d-r$   & $2d+r$ & $0 \leq r \leq d-1$  \\
\hline
 $d+2r+1$ &\ $2d-1-r$   & $3d+r$ &  $0 \leq r \leq d-2$  \\
\hline
\end{tabular}}
\caption{Langford sequence, $L_{d}^{2d-1}$.} \label{tablel}
\end{center}
\end{table}
\end{construction}
%\begin{lemma}\label{lldd}The pairs given by Construction~\ref{ldd} produce a Langford sequence with defect $d$ and order $2d-1$.\end{lemma}
%\\
Construction~\ref{ldd} is new and it is straightforward to check that Construction~\ref{ldd} gives a Langford sequence.
%%%%%%%%%%%%%%%%%%%%%%%%%%%%%%%%%%%%%%%    SKolem %%%%%%%%%%%%%%%%%%%%%%%%%%%%%%%%%%%%%%%%%%%%%%%%%%%%%%%%%

Notice that there is a different use of defect in Langford and near-Skolem sequences. In Langford sequences the defect ($d$) is the smallest integer in the sequence, but in near-Skolem sequences the defect ($m$) is the integer omitted from the sequence.
%A {\em hooked near-Skolem sequence} of order $n$ and defect $m$ is a sequence $m$-near $hnSm=(s_{1},s_{2},\ldots,s_{2n-1})$ is a hooked Skolem-like sequence of order $n-1$ with $H=\{1,2,\ldots, m-1,m+1,\ldots,n\}$.
%%%%%%%%%%%%%%%%%%%%%%%%%%%%%%%%%%%%%%%%%%%%%%%%%%%%% n=4m+1 %%%%%%%%%%%%%%%%%%%%%%%%%%%%%%%%%%%%%%%%%%%%%%%%%%%%%%%%%%%%%%%%%%%%%%%%%%%%%%%%%%%%%%%%%%%%%%%%%%%%%%%%%%%
\begin{construction}\label{ns11}From Table~\ref{ns1}, we can construct a hooked near-Skolem sequence of order $n=4m+1$ with $m \geq 3$ and defect $n-1$, (omitting row $5$ when $m=3$).\end{construction}
%\\ $\left(5,3,1,1,3,5,2,0,2\right).$
%\\ $\left(9,7,3,4,6,3,5,4,7,9,6,5,1,1,2,0,2\right).$\\
%Though they do not fit the construction of Table~\ref{ns1}, we note that for $n=5, \left(5,3,1,1,3,5,2,0,2\right)$ and $n=9, \left(9,7,3,4,6,3,5,4,7,9,$\newline$6,5,1,1,2,0,2$) are hooked near Skolem sequences of order $n=5$ and $n=9$, respectively, with defect $n-1$. Thus, hooked near Skolem sequences or order $4m+1$ with defect $4m$ exist for all $m$.
%%%%%%%%%%%%%%%%%%%%%%%%%%%%%%%%%%%%%%%%%%%%%%%%%%%%%%%%%%%%%%%%%%%%%%%%%%%%%%%%%%%%%%%%%%%%%%%%%%%%%%%%%%%%%%%%%%%%%%%%%%%%%%%%%%%%%%%%%%%%%%%%%%%%%%%%%%%%%%%%
%%%%%%%%%%%%%%%%%%%%%%%%%%%%%%%%%%%%%%%%%%%%%%%%%%%%% n=4m+3 %%%%%%%%%%%%%%%%%%%%%%%%%%%%%%%%%%%%%%%%%%%%%%%%%%%%%%%%%%%%%%%%%%%%%%%%%%%%%%%%%%%%%%%%%%%%%%%%%%%%%%%%%%%
%\begin{table}[H]\setlength\extrarowheight{6pt}\renewcommand{\arraystretch}{0.9}
%\begin{center}\setlength\extrarowheight{6pt}\renewcommand{\arraystretch}{0.9}
  %\begin{minipage}[b]{0.42\textwidth}
%\scalebox{0.85}{
\begin{table}[ht!]
\begin{center}%\setlength\extrarowheight{7pt}\renewcommand{\arraystretch}{0.9}
\scalebox{1.00}{
\begin{tabular}{|c|c|c|c|}
\hline
  $i$ & $a_i$ & $b_i$ &  \\
\hline
\hline
 $2r+1$ &\ $2m+1-r$   & $2m+2+r$ & $1 \leq r \leq 2m$  \\
\hline
 $4m-4$ &\ $2m+2$   & $6m-2$ &  $-$  \\
\hline
 $4m-2$ &\ $2m+1$   & $6m-1$ &  $-$  \\
\hline
 $2m+2r$ &\ $5m-r$   & $7m+r$ &  $0 \leq r \leq m-3$  \\
\hline
 $2r+4$ &\ $6m-r-3$   & $6m+r+1$ &  $0 \leq r \leq m-4$  \\
\hline
 $2m-2$ &\ $6m$   & $8m-2$ &  $-$  \\
\hline
 $1$ &\ $7m-2$   & $7m-1$ &  $-$  \\
\hline
 $2$ &\ $8m-1$   & $8m+1$ &  $-$  \\
\hline
\end{tabular}}
\caption{Hooked near-Skolem sequence from Construction~\ref{ns11}.}\label{ns1}
\end{center}
\end{table}
%\end{minipage}
%\quad
%\quad
%\begin{minipage}[b]{0.41\textwidth}
\begin{construction}\label{ns22}From Table~\ref{ns2}, we can construct a near-Skolem sequence of order $n=4m+3$ with $m \geq 2$ and defect $n-1$, (omitting row $4,5$ when $m=2$).\end{construction}
\begin{table}[ht!]
\begin{center}%\setlength\extrarowheight{7pt}\renewcommand{\arraystretch}{0.9}
\scalebox{1.00}{
\begin{tabular}{|c|c|c|c|}
\hline
  $i$ & $a_i$ & $b_i$ &  \\
\hline
\hline
 $2r+1$ &\ $2m+2-r$   & $2m+3+r$ & $1 \leq r \leq 2m+1$  \\
\hline
 $4m$ &\ $2m+2$   & $6m+2$ &  $-$  \\
\hline
 $4m-2$ &\ $2m+3$   & $6m+1$ &  $-$  \\
\hline
 $2m+2+2r$ &\ $5m+2-r$   & $7m+4+r$ &  $0 \leq r \leq m-3$  \\
\hline
 $2r+4$ &\ $6m-r$   & $6m+4+r$ &  $0 \leq r \leq m-3$  \\
\hline
 $1$ &\ $7m+2$   & $7m+3$ &  $-$  \\
\hline
 $2m$ &\ $6m+3$   & $8m+3$ &  $-$  \\
\hline
 $2$ &\ $8m+2$   & $8m+4$ &  $-$  \\
\hline
\end{tabular}}
\caption{Near-Skolem sequence from Construction~\ref{ns22}.}\label{ns2}
%\caption*{Near-Skolem sequence construction of order $n=4m+3$ with defect $n-1$.} \label{ns2}
%\end{minipage}
\end{center}
\end{table}
%\\ $\left(5,7,2,4,2,5,3,4,7,3,1,1\right).$\\
%\textbf{\textcolor[rgb]{1,0,0}{Note : I am unable to construct sequencexs for $n=5,7,9$.}}
%Though they do not fit the construction of Table~\ref{ns2}, we note that for $n=7, \left(5,7,2,4,2,5,3,4,7,3,1,1\right)$ is hooked near Skolem sequence with defect $n-1$. Thus, near Skolem sequences or order $4m+3$ with defect $4m+2$ exist for all $m$.
It is straightforward to check that Constructions~\ref{ns11} and~\ref{ns22} both produce the desired sequences. These constructions establish Lemma~\ref{hns}, which will be useful in Section~\ref{section3} and~\ref{section5} when labelling $C_{5}^{p}$ and $C_{3}^{t}C_{5}^{p}$.
\begin{lemma}\label{hns}If $n = 2k+1$, with $k \geq 5$, then there exists a (hooked) near-Skolem sequence of order $n,$ with defect $n-1$, which has no right endpoints in the positions $\left[1,k+2\right]$.\end{lemma}
%\begin{lemma}\label{hns1}If $n = 4m+1$, with $m \geq 3$, then exists a near hooked Skolem sequence of order $n,$ with defect $n-1$, which has no right endpoints in the positions $\left\{1,2,3,\ldots,2m+2\right\}$.\end{lemma}
%As a result of the above construction we obtained the following lemma. It is straightforward to check that the above construction it gives near Skolem sequence of order $n=4m+3$ with defect $n-1$.
%\begin{lemma}\label{hns2}If $n = 4m+3$, with $m \geq 2$, then exists a near hooked Skolem sequence of order $n,$ with defect $n-1$, which has no right endpoints in the positions $\left\{1,2,3,\ldots,2m+3\right\}$.\end{lemma}
%\\
%The necessary and sufficient conditions for the existence of Skolem, Langford and near Skolem sequences are given in \cite{crc}.\par
%It is straightforward to check that the constructions presented in Table~\ref{ns1} and Table~\ref{ns2}.
%%%%%%%%%%%%%%%%%%%%%%%%%%%%%%%%%%%%%%%%%%%%%%%%%%%%%%%%%%%%%%%%%%%%%%%%%%%%%%%%%%%%%%%%%%%%%%%%%%%%%%%%%%%%%%%%%%%%%%%%%%%%%%%%%%%%%%%%%%%%%%%%%%%%%%%%%%%%%%%%%%%%%%%%
%%%%%%%%%%%%%%%%%%%%%%%%%%%%%%%%%%%%%%%%%%%%%%% Two fold Skolem Like Sequence %%%%%%%%%%%%%%%%%%%%%%%%%%%%%%%%%%%%%%%%%%%%%%%%%%%%%%%%%%%%%%%%%%%%%%%%%%%%%%%%%%%%%%%%
%Here we extend the definition of Skolem-like sequences, we will define the two-fold Skolem-like sequences, these sequences will use each element four times.\\
We now extend our definition of Skolem-type sequences to sequences where every pair of elements in $H$ occurs exactly twice. A \textit{two-fold Skolem-type sequence} of order $n$ is a sequence $\left(s_{1},s_{2},\ldots,s_{4n}\right)$ of $4n$ positive integers such that the following conditions hold:
\begin{enumerate}
\item for a set $H$ of $n$ distinct positive integers, for every $p\in H,$  there exist exactly $2$ disjoint pairs $\left(s_{i},s_{j}\right)$ and $(s_{i'},s_{j'})$ such that $s_{i}=s_{j}=s_{i'}=s_{j'}=p$,
\item if $s_{i}=s_{j}=p$ and $s_{i'}=s_{j'}=p$ with $i < j$ and $i' < j'$, then $j-i=p$ and $j' - i'=p$.
%\item $s_{i},s_{j},s_{i'},s_{j'} \in K_{n}^{2}$ and $s_{i}=s_{j}=p$ with $s_{i'}=s_{j'}=p$.
\end{enumerate}
%%%%%%%%%%%%%%%%%%%%%%%%%%%%%%%%%%%%%%%%%%%%%%%%%%%%%%%%%%%%
%%%%%%%%%%%%%%%%%%%%%%%%%%%%%%%%%%%%%%%%%%%%%%%%%%%%%%%%%%%%%%%%%%%%%%%%%%%%%%%%%%%%%%%%%%%%%%%%%%%%%%%%%%%%%%%%%%%%%%%%%%%%%%%%%%%%%%%%%%%%%%%%%%%%%%%%%%%%%%%%%%%%%%%
For example, $\left(2,3,2,3,3,2,3,2\right)$ is a two-fold Skolem-type sequence of order $2$ with $H=\left\{2,3\right\}$.
%%%%%%%%%%%%%%%%%%%%%%%%%%%%%%%%%%%%%%%%%%%%%%%%%%%%%%%%%%%extended skolem like %%%%%%%%%%%%%%%%%%%%%%%%%%%%%%%%%%%%%%%%%%%%%%%%%%%%%%%%%%%%%%%%%%%%%%%%%%%%%%%%%%%%%%%%
%%%%%%%%%%%%%%%%%%%%%%%%%%%%%%%%%%%%%%%%%%%%%%%%%%%%%%%%%%%%%%%%%%%%%%%%%%%%%%%%%%%%%%%%%%%%%%%%%%%%%%%%%%%%%%%%%%%%%%%%%%%%%%%%%%%%%%%%%%%%%%%%%%%%%%%%%%%%%%%%%%%%%%%
\begin{construction}\label{like1}From Table~\ref{tablek}, we can construct a two-fold Skolem-type sequence construction of order $n \geq 1$, where $H=\left\{1\right\} \cup \left\{4i | 1 \leq i \leq n-1\right\}$.\end{construction}
% is a set of $n$ distinct positive integers.
%%%%%%%%%%%%%%%%%%%%%%%%%%%%%%%%%%%%%%%%%%%%%%%%%%%%%%%%%%%%%%%%%%%%%%%%%%%%%%%%%%%%%%%%%%%%%%%%%%%%%%%%%%%%%%%%%%%%%%%%%%%%%%%%%%%%%%%%%%%%%%%%%%%%%%%%%%%%%%%%%%%%%%%
%\begin{lemma}\label{like11}The pairs given by Construction~\ref{like1} produces a two-fold Skolem-like sequence.\end{lemma}
%%%%%%%%%%%%%%%%%%%%%%%%%%%%%%%%%%%%%%%%%%%%%%%%%%%%%%%%%%%%%%%%%%%%%%%%%%%%%%%%%%%%%%%%%%%%%%%%%%%%%%%%%%%%%%%%%%%%%%%%%%%%%%%%%%%%%%%%%%%%%%%%%%%%%%%%%%%%%%%%%%%%%%%
%\begin{table}[H]
%\begin{center}\setlength\extrarowheight{6pt}\renewcommand{\arraystretch}{1.4}
  %\begin{minipage}{0.4\textwidth}
%\scalebox{0.55}{
\begin{table}[ht!]
\begin{center}%\setlength\extrarowheight{7pt}\renewcommand{\arraystretch}{0.9}
\scalebox{1.00}{
\begin{tabular}{|c|c|c|c|c|}
\hline
  $j$ & $\left(c_{j},d_{j}\right)$ & $\left(e_{j},f_{j}\right)$ &  \\
\hline
\hline
 $1$ &\ $\left(2n-1,2n\right)$   & $\left(4n-1,4n\right)$ & $-$  \\
\hline
 $4r$ &\ $\left(2n-2r-1,2n+2r-1\right)$   & $\left(2n-2r,2n+2r\right)$ &  $1 \leq r \leq n-1$  \\
\hline
\end{tabular}}
\caption{Two-fold Skolem-type sequence from Construction~\ref{like1}.} \label{tablek}
\end{center}
\end{table}
For example, $\left(8,8,4,4,1,1,4,4,8,8,1,1\right)$ is a two-fold Skolem-type sequence of order $3$ with $H=\left\{1,4,8\right\}$, created using Construction~\ref{like1}.

A \textit{$\left(2n+1,2n+2\right)$-extended two-fold Skolem-type sequence} of order $n$ is a sequence $\left(\hat s_{1},\hat s_{2},\ldots,\hat s_{4n+2}\right)$ of $4n+2$ positive integers that satisfies the above conditions as well as $\hat s_{2n+1}=\hat s_{2n+2}=0$.

\begin{proposition}\label{like3}The following sequences are all two-fold Skolem-type sequences:

$C^{0}=\emptyset$, (the sequence of length $0$),

$C^{1}=\left(2,2,2,2\right)$, where $H=\left\{2\right\}$,

$C^{2}=\left(2,3,2,3,3,2,3,2\right)$, where $H=\left\{2,3\right\}$,

$C^{3}=\left(2,2,2,2,5,3,5,3,3,5,3,5\right)$, where $H=\left\{2,3,5\right\}$, and

$C^{4}=\left(6,6,2,2,2,2,6,6,5,3,5,3,3,5,3,5\right)$, where $H=\left\{2,3,5,6\right\}$.
\end{proposition}
%Every two-fold Skolem sequence, two-fold Langford sequence, and a double Langford sequence is a two-fold Skolem-like sequence.
By placing restrictions on the nature of $H$, we can define two-fold Skolem and Langford sequences, similar to the way we defined the Skolem and Langford sequences.

A \textit{double Skolem-type sequence} of order $l$ and defect $d$ is a sequence obtained by concatenating two existing Skolem-type sequences with the same defect and order, or by interlacing a hooked Skolem-type sequence with its reverse. A \textit{two-fold Skolem sequence} of order $n$ is a two-fold Skolem-type sequence with $H=\left[1,n\right]$. It is straightforward to show that two-fold Skolem sequences exist for any order $n$ by concatenating Skolem sequences (possibly interlacing their hooks). However, in Tables~\ref{table2fo} and~\ref{table2fe} we introduce a different construction for two-fold Skolem sequences which will be useful later.

\begin{table}[ht!]
\begin{center}%\setlength\extrarowheight{7pt}
\scalebox{1.00}{\renewcommand{\arraystretch}{1.1}
\begin{tabular}{|c|c|c|c|c|}
\hline
  $j$ & $\left(c_{j},d_{j}\right)$ & $\left(e_{j},f_{j}\right)$ &  \\
\hline
\hline
 $2r+1$ &\ $\left(\frac{n+1}{2}-r,\frac{n+3}{2}+r\right)$   & $\left(\frac{3n+3}{2}-r,\frac{3n+5}{2}+r\right)$ & $0 \!\leq\! r \!\leq\! \frac{n-1}{2}$  \\
\hline
 $n-1$ &\ $\left(2n+3,3n+2\right)$   & $\left(\frac{5n+5}{2},\frac{7n+3}{2}\right)$ &  $-$  \\
\hline
 $2r+2$ &\ $\left(\frac{5n+3}{2}-r,\frac{5n+3}{2}+r+2\right)$   & $\left(\frac{7n+1}{2}-r,\frac{7n+1}{2}+r+2\right)$ &  $0\! \leq \!r \!\leq\! \frac{n-5}{2}$  \\
\hline
%4& $2r+2$ &\ $-$   & $\left(\frac{7n+1}{2}-r,\frac{7n+1}{2}+r+2\right)$ &  $0 \leq r \leq \frac{n-5}{2}$  \\
%\hline
\end{tabular}}
\caption{Two-fold Skolem sequence construction of odd order $n$.}\label{table2fo}
\end{center}
\end{table}

\begin{table}[ht!]
\begin{center}%\setlength\extrarowheight{7pt}\renewcommand{\arraystretch}{1.1}
\scalebox{1.00}{\renewcommand{\arraystretch}{1.1}
\begin{tabular}{|c|c|c|c|c|}
\hline
  $j$ & $\left(c_{j},d_{j}\right)$ & $\left(e_{j},f_{j}\right)$ &  \\
\hline
\hline
 $2r+1$ &\ $\left(\frac{n}{2}-r,\frac{n+2}{2}+r\right)$   & $\left(\frac{3n}{2}-r,\frac{3n+2}{2}+r\right)$ & $0 \leq r \leq \frac{n-2}{2}$  \\
\hline
 $n$ &\ $\left(2n+1,3n+1\right)$   & $\left(\frac{5n+2}{2},\frac{7n+2}{2}\right)$ &  $-$  \\
\hline
$2r+2$ &\ $\left(\frac{5n}{2}-r,\frac{5n}{2}+r+2\right)$   & $\left(\frac{7n}{2}-r,\frac{7n}{2}+r+2\right)$ &  $0 \leq r \leq \frac{n-4}{2}$  \\
\hline
%4& $2r+2$ &\ $-$   & $\left(\frac{7n}{2}-r,\frac{7n}{2}+r+2\right)$ &  $0 \leq r \leq \frac{n-4}{2}$  \\
%\hline
\end{tabular}}
\caption{Two-fold Skolem sequence construction of even order $n$.}\label{table2fe}
\end{center}
\end{table}

A \textit{two-fold Langford sequence} of order $l$ and defect $d$ is a two-fold Skolem-type sequence of order $l$, $\left(l_{1},l_{2},\ldots,l_{4l}\right)$ with $H=\left[d,d+l-1\right]$. Again, it is straightforward to show that a two-fold Langford sequence can be obtained by concatenating (hooked) Langford sequences. In Construction~\ref{2foldlang}, we give a two-fold Langford sequence that is not constructed by concatenation that will be useful later.

\begin{construction}\label{2foldlang}From Table~\ref{table2l}, we can construct a two-fold Langford sequence with defect $6k-1$ and order $4k-1$, where $k \geq 1$.
\end{construction}
\begin{table}[ht!]
\begin{center}%\setlength\extrarowheight{7pt}\renewcommand{\arraystretch}{1.1}
\scalebox{1}{\renewcommand{\arraystretch}{1.1}
\begin{tabular}{|@{\hspace{.1cm}}c@{\hspace{.1cm}}|@{\hspace{.1cm}}c@{\hspace{.1cm}}|@{\hspace{.1cm}}c@{\hspace{.1cm}}|@{\hspace{.1cm}}c@{\hspace{.1cm}}|c}
\hline
  $j$ & $\left(c_{j},d_{j}\right)$ & $\left(e_{j},f_{j}\right)$ &  \\
\hline
\hline
 $10k\!-\!2\!-\!2r$ &\ $\left(r,10k\!-\!2\!-\!r\right)$   & $\left(2k\!-\!1\!+\!r,12k\!-\!3\!-\!r\right)$ & $1 \leq r \leq 2k\!-\!1$  \\
\hline
 $10k\!-\!1\!-\!2r$ &\ $\left(4k\!-\!2\!+\!r,14k\!-\!3\!-\!r\right)$   & $\left(6k\!-\!2\!+\!r,16k\!-\!3\!-\!r\right)$ &  $1 \leq r \leq 2k$  \\
\hline  
\end{tabular}}
\caption{Two-fold Langford sequence with defect $6k-1$ and order $4k-1$, where $k \geq 1$.} \label{table2l}
\end{center}
\end{table}

We summarize necessary and sufficient conditions for the existence of various useful Skolem-type sequences in Table~\ref{sbig} as an aid for the reader.
%\begin{tabular}{|>{\centering\arraybackslash}p{0.1\linewidth}|>{\centering\arraybackslash}p{1.3\linewidth}|>{\centering\arraybackslash}p{0.1\linewidth}|}
\begin{table}[ht!]
\large
\begin{center}\setlength\extrarowheight{2pt}\renewcommand{\arraystretch}{1.0}
\scalebox{.88}{
\begin{tabular}{|>{\raggedright\arraybackslash}p{0.25\linewidth}|>{\centering\arraybackslash}p{0.11\linewidth}|>{\raggedright\arraybackslash}p{0.56\linewidth}|>{\centering\arraybackslash}p{0.07\linewidth}|} % \left[   \right]
%\begin{tabular}{ | m{14cm} | m{7cm}| m{16cm}  }
\hline
Sequence & Order/ Defect & Necessary and Sufficient conditions & Ref.   \\ \hline
Skolem & $n$/-     &  $n\equiv0,1\ ({\rm mod}\ 4)$  & \cite{skolem}   \\ \hline
%5%%%%%%%%%%%%%%%%%%%%%%%%%%%%%%%%%%%%%%%%%%%%%%%%%%%%%%%%%%%%%%%%%%%%%%%%%%%%%%%%%%%%%%%%%%%%%%%%%%%%%%%%%%%%%%%%%%%%%%%%%%%%%%%%%%%%%%%%%%%%%%%%%%%%%%%%%%%%%%%%%%%%%
hooked Skolem & $n$/-     &  $n\equiv2,3\ ({\rm mod}\ 4)$  & \cite{okeefe}   \\ \hline	
%5%%%%%%%%%%%%%%%%%%%%%%%%%%%%%%%%%%%%%%%%%%%%%%%%%%%%%%%%%%%%%%%%%%%%%%%%%%%%%%%%%%%%%%%%%%%%%%%%%%%%%%%%%%%%%%%%%%%%%%%%%%%%%%%%%%%%%%%%%%%%%%%%%%%%%%%%%%%%%%%%%%%%%
Langford & $l$/$d$     & $l \geq 2d-1$,\newline $l \equiv 0,1 \pmod 4$ and $d$ is odd, or $l \equiv 0,3 \pmod 4$ and $d$ is even  & \cite{simpson}   \\ \hline
%5%%%%%%%%%%%%%%%%%%%%%%%%%%%%%%%%%%%%%%%%%%%%%%%%%%%%%%%%%%%%%%%%%%%%%%%%%%%%%%%%%%%%%%%%%%%%%%%%%%%%%%%%%%%%%%%%%%%%%%%%%%%%%%%%%%%%%%%%%%%%%%%%%%%%%%%%%%%%%%%%%%%%%	
hooked Langford & $l$/$d$     & $l\left(l-2d+1\right)+2 \geq 0$,\newline $l \equiv 2,3 \pmod 4$ and $d$ is odd, or $l \equiv 1,2 \pmod 4$ and $d$ is even  & \cite{simpson}   \\ \hline	
%5%%%%%%%%%%%%%%%%%%%%%%%%%%%%%%%%%%%%%%%%%%%%%%%%%%%%%%%%%%%%%%%%%%%%%%%%%%%%%%%%%%%%%%%%%%%%%%%%%%%%%%%%%%%%%%%%%%%%%%%%%%%%%%%%%%%%%%%%%%%%%%%%%%%%%%%%%%%%%%%%%%%%%	
$m$-near-Skolem & $n$/-     &$n \equiv 0,1 \pmod 4$ and $m$ is odd, or $n \equiv 2,3 \pmod 4$ and $m$ is even  & \cite{shalaby}   \\ \hline	
%5%%%%%%%%%%%%%%%%%%%%%%%%%%%%%%%%%%%%%%%%%%%%%%%%%%%%%%%%%%%%%%%%%%%%%%%%%%%%%%%%%%%%%%%%%%%%%%%%%%%%%%%%%%%%%%%%%%%%%%%%%%%%%%%%%%%%%%%%%%%%%%%%%%%%%%%%%%%%%%%%%%%%%	
hooked $m$-near-Skolem& $n$/-     &$n \equiv 0,1 \pmod 4$ and $m$ is even, or $n \equiv 2,3 \pmod 4$ and $m$ is odd  & \cite{shalaby}   \\ \hline	
%5%%%%%%%%%%%%%%%%%%%%%%%%%%%%%%%%%%%%%%%%%%%%%%%%%%%%%%%%%%%%%%%%%%%%%%%%%%%%%%%%%%%%%%%%%%%%%%%%%%%%%%%%%%%%%%%%%%%%%%%%%%%%%%%%%%%%%%%%%%%%%%%%%%%%%%%%%%%%%%%%%%%%%	
$m$-fold Skolem & $n$/-     &$n \equiv 0,1 \pmod 4$ and any $m$, or $n \equiv 2,3 \pmod 4$ and $m$ is even  & \cite{baker}   \\ \hline
%5%%%%%%%%%%%%%%%%%%%%%%%%%%%%%%%%%%%%%%%%%%%%%%%%%%%%%%%%%%%%%%%%%%%%%%%%%%%%%%%%%%%%%%%%%%%%%%%%%%%%%%%%%%%%%%%%%%%%%%%%%%%%%%%%%%%%%%%%%%%%%%%%%%%%%%%%%%%%%%%%%%%%%	
hooked $m$-fold Skolem & $n$/-     &$n \equiv 2,3 \pmod 4$ and $m$ is odd  & \cite{baker}   \\ \hline		
\end{tabular}}
\caption{Summary of necessary and sufficient conditions for the existence of various Skolem-type sequences.}\label{sbig}
\end{center}
\end{table}

\section{Labellings from Skolem-type sequences}\label{section3}
In this section, we use Skolem-type sequences to label variable windmills of different orders.
%%%%%%%%%%%%%%%%%%%%%%%%%%%%%%%%%%%%%%%%%%%%%%%%%%%%%%%%%%%%%%%%%%%%%%%%%%%%%%%%%%%%%%%%%%%%%%%%%%%%%%%%%%%%%%%%%%%%%%%%%%%%%%%%%%%%%%%%%%%%%%%%%%%%%%%%%%%%%%%%%%%%%%%
%%%%%%%%%%%%%%%%%%%%%%%%%%%%%%%%%%%%%%%%%%%%%%%%%%%%% Constructions for 3-4 Windmills t >= s %%%%%%%%%%%%%%%%%%%%%%%%%%%%%%%%%%%%%%%%%%%%%%%%%%%%%%%%%%%%%%%%%%%%%%%%%%
%%%%%%%%%%%%%%%%%%%%%%%%%%%%%%%%%%%%%%%%%%%%%%%%%%%%%%%%%%%%%%%%%%%%%%%%%%%%%%%%%%%%%%%%%%%%%%%%%%%%%%%%%%%%%%%%%%%%%%%%%%%%%%%%%%%%%%%%%%%%%%%%%%%%%%%%%%%%%%%%%%%%%%%
\subsection{Constructions using Skolem-type sequences}
To begin, we discuss how to label $C_{3}^{t}$ and $C_{5}^{p}$ by using Skolem-type sequences. In \cite{yang}, Yang et~al.~proved that the $C_{5}^{p}$ is graceful when $p \equiv 0,3\  ({\rm mod}\ 4)$. Bermond in \cite{Bermond} proved that the $C_{3}^{t}$ is graceful when $t \equiv 0,1\  ({\rm mod}\ 4)$.

\begin{construction}\label{skolem01}
From a Skolem-type sequence or a hooked Skolem-type sequence of order $t$, construct the pairs $(a_i,b_i)$ such that $b_i-a_i=i$ for $1\leq i\leq t$. From these pairs, along with an arbitrary positive integer $c$, form one of the following sets of triples
\begin{enumerate}
\item {\label{skolem011}} $\left(0,a_i+c,b_i+c\right)$, $1\leq i\leq t$, or
\item {\label{skolem012}} $\left(0,i,b_i+c\right)$, $1\leq i\leq t$.
\end{enumerate}
The triples give a labelling of a $C_{3}^{t}$, with the central vertex labelled $0$.\end{construction}

\begin{lemma}\label{l1}Let $c \geq t$ be an arbitrary positive integer.
\hspace*{1em}
\begin{enumerate}
	\item Using any Skolem sequence of order $t$, with Construction~\ref{skolem01}(\ref{skolem011}) gives the edge labels $\left[1,t\right] \cup \left[c+1,c+2t\right]$, and the vertex labels from $\left\{0\right\} \cup \left[c+1,c+2t\right]$, such that each nonzero label occurs exactly once, and
	\item using any hooked Skolem sequence of order $t$, with Construction~\ref{skolem01}(\ref{skolem011}) gives the edge labels $\left[1,t\right] \cup \left[c+1,c+2t-1\right] \cup \left\{c+2t+1\right\}$, and the vertex labels $\left\{0\right\} \cup \left[c+1,c+2t-1\right]  \cup \left\{c+2t+1\right\}$, such that each nonzero label occurring exactly once, and
	\item  using any Skolem sequence with Construction~\ref{skolem01}(\ref{skolem012}) gives the edge labels $\left[1,t\right] \cup \left[c+1,c+2t\right]$ and the vertex labels from $\left\{0\right\} \cup \left[1,t\right] \cup \left[c+2,c+2t\right]$, each nonzero label occurring exactly once, and
\item using any hooked Skolem sequence with Construction~\ref{skolem01}(\ref{skolem012}) gives the edge labels $\left[1,t\right] \cup \left[c+1,c+2t-1\right] \cup \left\{c+2t+1\right\}$ and the vertex labels from $\left\{0\right\} \cup \left[1,t\right] \cup \left[c+2,c+2t-1,\right]  \cup \left\{c+2t+1\right\}$, each nonzero label occurring exactly once.
\end{enumerate}
\end{lemma}
\begin{pf}Start with a Skolem sequence $S_t$ of order $t$, and construct the triples $\left(0, a_i+c, b_i+c\right)_{i=1}^t$ by using Construction~\ref{skolem01}(\ref{skolem011}). Since $1 \leq a_i \leq 2t-1$ and $2 \leq b_i \leq 2t$, and all $a_i$ and $b_i$ are distinct, the vertex labels are $\left\{0\right\} \cup \left[c+1,c+2t\right]$, where no vertex label other than $0$ is repeated. Considering the edge labels, we see that they are $b_{i}-a_{i} = i$, $a_{i}+c$, and $b_{i}+c$. Since by construction, $1\leq i\leq t$ and all the $a_{i}$ and $b_{i}$ are distinct, we obtain edge labels $\left[1,t\right] \cup \left[c+1,c+2t\right]$, a union of disjoint sets since $c \geq t$, all of which are distinct.

A similar argument holds for hooked Skolem sequences, and for each case of Construction~\ref{skolem01}(\ref{skolem012}).\end{pf}
We make note of one more special case of a Skolem-type sequence used with Construction~\ref{skolem01}(\ref{skolem011}).

\begin{lemma}\label{llang1}Using a Langford sequence with Construction~\ref{skolem01}(\ref{skolem011}) and $c \geq d+l-1$ gives edge labels $\left[d,d+l-1\right] \cup \left[c+1,c+2l\right]$ and vertex labels from $\left\{0\right\} \cup \left[c+1,c+2l\right]$, each nonzero label occurring exactly once.\end{lemma}

To illustrate Constructions~\ref{skolem01}(\ref{skolem011}) and~\ref{skolem01}(\ref{skolem012}), consider for example, $(3$, $1$, $1$, $3$, $2$, $0$, $2)$, a hooked Skolem sequence of order $3$ which yields the pairs $(2,3)$, $(5,7)$, $(1,4)$. Letting $c=3$, these pairs yield the triples to near gracefully label a $C_{3}^{3}$:
$\left(0,5,6\right),\left(0,8,10\right)$, and $\left(0,4,7\right)$ by Construction~\ref{skolem01}(\ref{skolem011}), and
$\left(0,1,6\right),\left(0,2,10\right)$, and $\left(0,3,7\right)$ by Construction~\ref{skolem01}(\ref{skolem012}).
Taking $c=t$ and considering only Skolem and hooked Skolem sequences in this construction, we can use the resulting triples to (near) gracefully label Dutch windmills ($C_{3}^{t}$) for any $t \geq 1$, as well as some related graphs. See \cite{dyer} and \cite{ahmad}.

Koh, Rogers, Lee, and Toh \cite{koh} conjectured that $C_{n}^{t}$ is graceful if and only if $nt \equiv 0,3\  ({\rm mod}\ 4)$. In 2005, Yang et al. in \cite{yang} have shown the conjecture true for $n = 5.$ In Theorem~\ref{thirdthm}, we prove that $G=C_{5}^{p}$ is near graceful when $p\equiv1,2\ ({\rm mod}\ 4)$, and verify that $G$ is graceful when $p \equiv 0,3\  ({\rm mod}\ 4)$ through the use of (hooked) Skolem and (hooked) near-Skolem sequences. Then, in Section 5, we use this construction to prove (near) graceful labellings exist for families of $C_{3}^{t}C_{5}^{p}$. In Construction~\ref{c5a} we will use (hooked) Skolem and (hooked) near-Skolem sequences together to form the $5$-tuples to label a $C_{5}^{p}$.% In Theorem~\ref{thirdthm}, we verify that $G=C_{5}^{p}$ is near graceful if and only if $p\equiv1,2\ ({\rm mod}\ 4)$ and $G$ is graceful when $p \equiv 0,3\  ({\rm mod}\ 4)$.
\begin{construction}\label{c5a}
Let $p > 0$. Given a (hooked) Skolem sequence $S_{1}$ of order $p=4m+k$ where $0 \leq k \leq 3$, construct the pairs $(a_i,b_i)$ such that $b_i-a_i=i$ for $1\leq i\leq p$. Select Skolem sequence $S_{2}$ and construct the pairs $(c_j,d_j)$ as follows.
\begin{enumerate}
	\item For $k=0$, select the Skolem sequence of order $2p$ using Table~\ref{table1}, in Appendix A.
	
	\item For $k=1$, select the hooked Skolem sequence of order $2p$ using Table~\ref{table3}, in Appendix A.
	
	\item For $k=2$, select the hooked near-Skolem sequence of order $2p+1$ using Table~\ref{ns1}.
	
	\item For $k=3$, select the near-Skolem sequence of order $2p+1$ using Table~\ref{ns2}.
\end{enumerate}
Then form the $5$-tuple $\left(0,d_{b_i}+p,b_i,a_i,d_{a_i}+p\right)$ for each $1 \leq i \leq p$.
\end{construction}
Note that by construction, for $k=0$ and $k=1$, no right endpoints occur in positions $1$ to $p$ of $S_2$. Similarly, $S_2$ has no right endpoints in positions $1,2,\ldots,p-1$, nor in $p+1$ when $k=2$ and $k=3$.
\begin{theorem}\label{thirdthm} If $G=C_{5}^{p}$, then $G$ is near graceful when $p \equiv 1,2\ ({\rm mod}\ 4)$ and $G$ is graceful when $p \equiv 0,3\ ({\rm mod}\ 4)$.\end{theorem}
\begin{pf} \\
\textbf{Case 1:} Let $p \equiv 1\ ({\rm mod}\ 4)$. Form the $5$-tuples $\left(0,d_{b_i}+p,b_i,a_i,d_{a_i}+p\right)$ for each $1 \leq i \leq p$ as indicated in Construction~\ref{c5a}. We begin by considering the vertex labels used by the $5$-tuples.

% Let $p \equiv 1\ ({\rm mod}\ 4)$. Start with a Skolem sequence of order $p=4m+1$ and construct the pairs $\left(a_i,b_i\right)$ with $1 \leq i \leq p$ as indicated in Construction~\ref{c5a}. Also, form a hooked Skolem sequence of order $2p$ from Table~\ref{table3} and construct the pairs $\left(c_j,d_j\right)$ with $1 \leq j \leq 2p$ as indicated in Construction~\ref{c5a}. Form the $5$-tuples $\left(0,\left(d_{b_i}\right)+p,b_i,a_i,\left(d_{a_i}\right)+p\right)$ for each $1 \leq i \leq p$.
Note that from the Skolem sequence $S_{1}$ the entries $a_i$ and $b_i$ in the third and fourth entries of the $5$-tuples give the distinct numbers $\left[1,2p\right]$. As we know, there are no right endpoints $\left(d_j\right)$ in the first $p$ positions in the hooked Skolem sequence. The first right endpoint is in position $p+1$, so the set of possible positions for right endpoints are $\left[p+1,4p-1\right]  \cup \left\{4p+1\right\}.$ Thus, the second and fifth entries of the $5$-tuples are all elements of $\left[2p+1,5p-1\right]  \cup \left\{5p+1\right\}.$ In Skolem sequences, $a_i$ and $b_i$ are distinct, so in the hooked Skolem sequence, $d_{a_i}$ and $d_{b_i}$ are distinct. Further, we know that the entries $a_i$ and $b_i$ on the third and fourth entries of the $5$-tuples give the distinct numbers $\left[1,2p\right]$, and $d_{a_i},d_{b_i} \geq p+1$ so the minimum value of the second and fifth entries of the $5$-tuples is $2p+1$. Therefore, all the nonzero entries of the $5$-tuples are distinct.
%The Skolem sequence $S_{2}$ contains no right endpoints in the first $p$ positions, with the first right end point occuring in position $p+1$. The union of $\left\{d_{a_i},d_{b_i}\right\}_{i=1}^{p}$ is the set of all right endpoints, and so a subset of $\left[p+1,4p\right].$ From the $5$-tuples, we use the right endpoints as indicated in the second and fifth entries. So we have to add $p$ to each element based on our $5$-tuples so we obtain a subset of $\left[2p+1,5p-1\right]  \cup \left\{5p+1\right\}$ and since $a_i$ and $b_i$ subset of $\left[1,2p\right]$ then we all the entries of the $5$-tuples are distinct.
% so we obtain the following distinct elements $\left\{p+1,p+2,\ldots,4p-1,4p+1\right\}.$ From the $5$-tuples we use the right endpoints as indicated in the second and fifth entries. We have to add $p$ for each element so we obtain $\left\{2p+1,2p+2,\ldots,5p-1,5p+1\right\}.$ In a Skolem sequence, the $a_i,b_i,c_i$, and $d_i$ are all unique, and the only repeated element over these $5$-tuples is $0$.

From the above discussion, it is clear that the only vertex label repeated is $0$ ($p$ times), and that all vertices are distinct and come from the set $\left[0,5p-1\right]  \cup \left\{5p+1\right\}$.

%%%%%%%%%%%%%%%%%%%%%%%%%%%%%%%%%%%%%%%%%%%%%%%%%%%%%%%%%%%%%%%%%%%%%%%%%%%%%%%%%%%%%%%%%%%%%%%%%%%%%%%%%%%%%%%%%%%%%%%%%%%%%%%%%%%%%%%%%%%%%%%%%%%%%%%%%%%%%%%%
%%%%%%%%%%%%%%%%%%%%%%%cccccccccccccccccccccccccccccccccccccccccccccccccccccccccccccccccccccccccccccccccccccccccccccccccccccccccccccccccccccccccccccccccccc
We now consider the edge labels defined by the difference between subsequent entries (taken cyclically) in the $5$-tuple $\left(0,d_{b_i}+p,b_i,a_i,d_{a_i}+p\right)$.

Since $b_i - a_i = i$, these differences are all distinct and comprise the set $\left[1,p\right]$. Based on our previous discussion, the differences between $d_{b_i}+p$ and $0$ and the differences between $d_{a_i}+p$ and $0$ give distinct numbers from the set $\left[2p+1,5p-1\right]  \cup \left\{5p+1\right\}$. Considering the remaining differences, we see that $\left(d_{b_i} + p\right) - b_i= \left(d_{b_i} - b_i\right) +p = c_{b_i} + p$ and $\left(d_{a_i} + p\right) - a_i= \left(d_{a_i} - a_i\right) +p = c_{a_i} + p$. These differences are all distinct numbers in the set $\left[p+1,5p-1\right]$. % Since$\cup \left\{a_i,b_i\right\}=\left\{1,2,\ldots,?\right\}$
 Since $\bigcup\limits_{i=1}^{p} \left\{a_i,b_i\right\} = \left[1,p\right]$, then $\bigcup\limits_{i=1}^{p} \left\{c_{a_i},c_{b_i},d_{a_i},d_{b_i}\right\} = \left[1,4p-1\right]  \cup \left\{4p+1\right\}$, and hence $\bigcup\limits_{i=1}^{p} \left\{c_{a_i}+p,c_{b_i}+p,d_{a_i}+p,d_{b_i}+p\right\}$ is exactly $\left[p+1,5p-1\right] \cup \left\{5p+1\right\}$, since all the $c_j$ and $d_j$ are distinct.
%In a Skolem sequence, the $a_i,b_i,c_i$ and $d_i$ are all unique, the only repeated element over these $5$-tuples is $0$. If we take the union of the above sets we get the edge labels $\left\{1,2,\ldots,5p-1,5p+1\right\}$.$\bigcup^{p}_{i=1}$$\bigcup\limits_{i=1}^{p} \left\{a,b\right\}$

From the above discussion, it is clear that all edges are distinct and are exactly the set $\left[1,5p-1\right]  \cup \left\{5p+1\right\}$. We can conclude that since the vertex labels are a subset of $\left[0,5p-1\right]  \cup \left\{5p+1\right\}$, each nonzero label occurring exactly once, and the edge labels are exactly $\left[1,5p-1\right]  \cup \left\{5p+1\right\}$, $C_{5}^{p}$ can be near gracefully labelled when $p \equiv 1\ ({\rm mod}\ 4)$.
%%%%%%%%%%%%%%%%%%%%%%%%%%%%%%%%%%%%%%%%%%%%%%%%%%%%%%%%%%%%%%%%%%%%%%%%%%%%%%%%%%%%%%%%%%%%%%%%%%%%%%%%%%%%%%%%%%%%%%%%%%%%%%
%We will prove one case only. The other cases follow similarly.

\textbf{Case 2:} The case $p \equiv 0\ ({\rm mod}\ 4)$, is proved similarly to Case 1. Use a Skolem sequence to construct the pairs $a_i$ and $b_i$. Instead of a hooked Skolem sequence, use a Skolem sequence to construct the pairs $d_{a_i}$ and $d_{b_i}$. The resulting $C_{5}^{p}$ can be gracefully labelled when $p \equiv 0\ ({\rm mod}\ 4)$.
%%%%%%%%%%%%%%%%%%%%%%%%%%%%%%%%%%%%%%%%%%%%%%%%%%%%%%%%%%%%%%%%%%%%%%%%%%%%%%%%%%%%%%%%%%%%%%%%%%%%%%%%%%%%%%%%%%%%%%%%%%%%%%

\textbf{Case 3:} If $p \equiv 2\ ({\rm mod}\ 4)$  and $p > 2$, the statement is proved similarly to Case 1. Instead of a Skolem sequence, use a hooked Skolem sequence to construct the pairs $a_i$ and $b_i$. Instead of a hooked Skolem sequence, use a hooked near-Skolem sequence to construct the pairs $d_{a_i}$ and $d_{b_i}$. We know such a sequence exists by Lemma~\ref{hns}. For $p=2$, use the $5$-tuples $\left(0,11,2,9,1\right)$ and $\left(0,6,3,7,5\right)$. Therefore, $C_{5}^{p}$ can be near gracefully labelled when $p \equiv 2\ ({\rm mod}\ 4)$.

%%%%%%%%%%%%%%%%%%%%%%%%%%%%%%%%%%%%%%%%%%%%%%%%%%%%%%%%%%%%%%%%%%%%%%%%%%%%%%%%%%%%%%%%%%%%%%%%%%%%%%%%%%%%%%%%%%%%%%%%%%%%%%
\textbf{Case 4:} If $p \equiv 3\ ({\rm mod}\ 4)$ and $p > 3$, the statement is proved similarly to Case 1. Instead of a Skolem sequence, use a hooked Skolem sequence to construct the pairs $a_i$ and $b_i$. Instead of a hooked Skolem sequence, use a near-Skolem sequence to construct  the pairs $d_{a_i}$ and $d_{b_i}$. We know such a sequence exists by Lemma~\ref{hns}. For $p=3$, use the $5$-tuples $\left(0,15,1,14,12\right)$, $\left(0,5,6,3,10\right)$ and $\left(0,9,13,2,8\right)$. Therefore, $C_{5}^{p}$ can be gracefully labelled when $p \equiv 3\ ({\rm mod}\ 4)$.
\end{pf}
\subsection{Two-fold Skolem-type sequence constructions}\label{subsec2}
In this section we discuss how to label $C_{4}^{s}$ using two-fold Skolem-type sequences.
%are the entries where $j$ occurs in the two-fold Skolem sequence,     and $t$ is a fixed positive integer
%%%%%%%%%%%%%%%%%%%%%%%%%%%%%%%%%%%%%%%%%%%%%%%%%%%%%%%%%%%%%%%%%%%%%%%%%%%%%%%%%%%%%%%%%%%%%%%%%%%%%%%
\begin{construction}\label{twofoldskolem03}
From a two-fold Skolem-type sequence of order $s$, construct the pairs of the form $\left(c_{j},d_{j}\right)$ and $\left(e_{j},f_{j}\right)$ where $c_{j},d_{j},e_{j}$, and $f_{j}$ are the entries where $j$ occurs in the two-fold Skolem-type sequence, $c_{j} < d_{j}, e_{j} < f_{j},$ with $j\in H$ and $d_{j}-c_{j}=f_{j}-e_{j}=j$. From these pairs we can obtain $s$ quadruples of the form $\left(0,d_{j}+c,j,f_{j}+c\right)$, where $c$ is a fixed positive integer. These quadruples admit a labelling of $C_{4}^{s}$.\end{construction}
In this paper, we will use Construction~\ref{twofoldskolem03} with a variety of two-fold Skolem-type sequences which will produce edge and vertex labels as given in Table~\ref{2big}. The symbols $t$ and $c$ are constants, $s$ and $l$ are the order of the sequences, and $d$ is the defect of the Langford sequence.
\begin{table}[ht!]
%\large
\begin{center}\setlength\extrarowheight{8pt}\renewcommand{\arraystretch}{1.0}
\scalebox{0.80}{
\begin{tabular}{|>{\raggedright\arraybackslash}p{0.48\linewidth}|>{\centering\arraybackslash}p{0.26\linewidth}|>{\raggedright\arraybackslash}p{0.4\linewidth}|} % \left[   \right]
%begin{tabular}{ | m{12cm} | m{7cm}| m{10cm}| }
\cline{1-3}
Sequence & Edge labels are: & Vertex labels from:   \\ \cline{1-3}
Two-fold Skolem sequences of order $s \leq c+1$     &    $\left[c+1,c+4s\right]$         &  $\left[0,s\right] \cup \left[c+2,c+4s\right]$               \\ \cline{1-3}
%5%%%%%%%%%%%%%%%%%%%%%%%%%%%%%%%%%%%%%%%%%%%%%%%%%%%%%%%%%%%%%%%%%%%%%%%%%%%%%%%%%%%%%%%%%%%%%%%%%%%%%%%%%%%%%%%%%%%%%%%%%%%%%%%%%%%%%%%%%%%%%%%%%%%%%%%%%%%%%%%%%%%%%	
Two-fold Skolem sequences from Table~\ref{table2fo} of order $s \leq 2c+1$ and $s$ is odd       &$\left[c+1,c+4s\right]$             & $\left[0,s\right]  \cup \left[c+\frac{s+3}{2},c+4s\right]$               \\ \cline{1-3}
%5%%%%%%%%%%%%%%%%%%%%%%%%%%%%%%%%%%%%%%%%%%%%%%%%%%%%%%%%%%%%%%%%%%%%%%%%%%%%%%%%%%%%%%%%%%%%%%%%%%%%%%%%%%%%%%%%%%%%%%%%%%%%%%%%%%%%%%%%%%%%%%%%%%%%%%%%%%%%%%%%%%%%%	
Two-fold Skolem sequences from Table~\ref{table2fe} of order $s \leq 2c$ and $s$ is even        &  $\left[c+1,c+4s\right]$           & $\left[0,s\right]  \cup \left[c+\frac{s+2}{2},c+4s\right]$               \\ \cline{1-3}
				%5%%%%%%%%%%%%%%%%%%%%%%%%%%%%%%%%%%%%%%%%%%%%%%%%%%%%%%%%%%%%%%%%%%%%%%%%%%%%%%%%%%%%%%%%%%%%%%%%%%%%%%%%%%%%%%%%%%%%%%%%%%%%%%%%%%%%%%%%%%%%%%%%%%%%%%%%%%%%%%%%%%%%%	
Double Langford sequence from Table~\ref{tablel} with defect $d \leq c+1$ and order $l$      & $\left[c+1,c+4l\right]$            &$\left\{0\right\} \cup \left[d,3d-2\right] \cup \left[c+2d,c+4l\right]$                \\ \cline{1-3}
				%5%%%%%%%%%%%%%%%%%%%%%%%%%%%%%%%%%%%%%%%%%%%%%%%%%%%%%%%%%%%%%%%%%%%%%%%%%%%%%%%%%%%%%%%%%%%%%%%%%%%%%%%%%%%%%%%%%%%%%%%%%%%%%%%%%%%%%%%%%%%%%%%%%%%%%%%%%%%%%%%%%%%%%	

Double Skolem sequences of order $s \leq 2c+2$ from Tables~\ref{table1}, \ref{table3}         & $\left[c+1,c+4s\right]$            & $\left[0,s\right] \cup \left[c+\frac{s+4}{2},c+4s\right]$               \\ \cline{1-3}
				%5%%%%%%%%%%%%%%%%%%%%%%%%%%%%%%%%%%%%%%%%%%%%%%%%%%%%%%%%%%%%%%%%%%%%%%%%%%%%%%%%%%%%%%%%%%%%%%%%%%%%%%%%%%%%%%%%%%%%%%%%%%%%%%%%%%%%%%%%%%%%%%%%%%%%%%%%%%%%%%%%%%%%%	

Double Skolem sequences of order $s \leq 2c+1$ from Tables~\ref{table2}, \ref{table4}         & $\left[c+1,c+4s\right]$            & $\left[0,s\right] \cup \left[c+\frac{s+3}{2},c+4s\right]$               \\ \cline{1-3}
				%5%%%%%%%%%%%%%%%%%%%%%%%%%%%%%%%%%%%%%%%%%%%%%%%%%%%%%%%%%%%%%%%%%%%%%%%%%%%%%%%%%%%%%%%%%%%%%%%%%%%%%%%%%%%%%%%%%%%%%%%%%%%%%%%%%%%%%%%%%%%%%%%%%%%%%%%%%%%%%%%%%%%%%	

Two-fold Skolem-type sequence of order $s \leq \frac{(c+4)}{2}$ from Table~\ref{tablek}       & $\left[c+1,c+4s\right]$            &$\left\{0\right\} \cup \left[c+2s+1,c+4s\right] \cup H$                \\ \cline{1-3}
				%5%%%%%%%%%%%%%%%%%%%%%%%%%%%%%%%%%%%%%%%%%%%%%%%%%%%%%%%%%%%%%%%%%%%%%%%%%%%%%%%%%%%%%%%%%%%%%%%%%%%%%%%%%%%%%%%%%%%%%%%%%%%%%%%%%%%%%%%%%%%%%%%%%%%%%%%%%%%%%%%%%%%%%	

%Two-fold Skolem-like sequence from Table~\ref{tablekt} of order $s$         &$\left[t+c+1,t+c+4s+2\right] \setminus$    $\left\{t+c+2s+1,t+c+2s+2\right\}$ &  $\left\{0\right\} \cup \left[t+c+2s+3,t+c+4s+2\right] \cup H$               &  \\ \cline{1-3}
				%5%%%%%%%%%%%%%%%%%%%%%%%%%%%%%%%%%%%%%%%%%%%%%%%%%%%%%%%%%%%%%%%%%%%%%%%%%%%%%%%%%%%%%%%%%%%%%%%%%%%%%%%%%%%%%%%%%%%%%%%%%%%%%%%%%%%%%%%%%%%%%%%%%%%%%%%%%%%%%%%%%%%%%	

Two-fold Skolem-type sequence $C^{1}$, with $c \geq 0$    &$\left[c+1,c+4\right]$             &$\left\{0,2\right\} \cup \left[c+3,c+4\right]$                \\ \cline{1-3}
				%5%%%%%%%%%%%%%%%%%%%%%%%%%%%%%%%%%%%%%%%%%%%%%%%%%%%%%%%%%%%%%%%%%%%%%%%%%%%%%%%%%%%%%%%%%%%%%%%%%%%%%%%%%%%%%%%%%%%%%%%%%%%%%%%%%%%%%%%%%%%%%%%%%%%%%%%%%%%%%%%%%%%%%	

Two-fold Skolem-type sequence $C^{2}$, with $c \geq 1$ &$\left[c+1,c+8\right]$ &$\left\{0,2,3\right\} \cup \left\{c+3,c+5,c+7,c+8\right\}$   \\ \cline{1-3}
				%5%%%%%%%%%%%%%%%%%%%%%%%%%%%%%%%%%%%%%%%%%%%%%%%%%%%%%%%%%%%%%%%%%%%%%%%%%%%%%%%%%%%%%%%%%%%%%%%%%%%%%%%%%%%%%%%%%%%%%%%%%%%%%%%%%%%%%%%%%%%%%%%%%%%%%%%%%%%%%%%%%%%%%	

Two-fold Skolem-type sequence $C^{3}$, with $c \geq 3$     & $\left[c+1,c+12\right]$             &$\left\{0,2,3,5\right\} \cup \left[c+3,c+4\right]$\newline $\cup \left[c+9,c+12\right]$               \\ \cline{1-3}
		%5%%%%%%%%%%%%%%%%%%%%%%%%%%%%%%%%%%%%%%%%%%%%%%%%%%%%%%%%%%%%%%%%%%%%%%%%%%%%%%%%%%%%%%%%%%%%%%%%%%%%%%%%%%%%%%%%%%%%%%%%%%%%%%%%%%%%%%%%%%%%%%%%%%%%%%%%%%%%%%%%%%%%%	

Two-fold Skolem-type sequence $C^{4}$, with $c \geq 2$     & $\left[c+1,c+16\right]$             &$\left\{0,2,3,5,6\right\} \cup \left[c+5,c+8\right] $\newline $\cup \left[c+13,c+16\right]$                \\ \cline{1-3}
		%5%%%%%%%%%%%%%%%%%%%%%%%%%%%%%%%%%%%%%%%%%%%%%%%%%%%%%%%%%%%%%%%%%%%%%%%%%%%%%%%%%%%%%%%%%%%%%%%%%%%%%%%%%%%%%%%%%%%%%%%%%%%%%%%%%%%%%%%%%%%%%%%%%%%%%%%%%%%%%%%%%%%%%	

Two-fold Langford sequences from Table~\ref{table2l} with $d=6k-1$ and $l=4k-1$, and $c \geq 2k-1$     & $\left[c+1,c+16k-4\right]$             &$\left\{0\right\}\cup \left[6k-1,10k-3\right]$\newline $\cup \left[c+8k-1,c+16k-4\right]$               \\ \cline{1-3}

\end{tabular}}
\caption{Summary of results of Construction~\ref{twofoldskolem03} with a variety of two-fold Skolem-type sequences.}\label{2big}
\end{center}
\end{table}
\begin{lemma}\label{22big}The results in Table~\ref{2big} are correct, with each nonzero label occurring exactly once.\end{lemma}
\begin{pf}We will prove that the result in the first row of Table~\ref{2big} is correct. The other rows follow similarly.

From a two-fold Skolem sequence $S_{s}^{2}$ made up of pairs $\left(c_{j},d_{j}\right)$ and $\left(e_{j},f_{j}\right)$, construct the quadruples $\left(0,d_{j}+c,j,f_{j}+c\right)_{j=1}^s$ as given by Construction~\ref{twofoldskolem03} with $c \geq s-1$. As we know $2 \leq d_j \leq 4s-1$ and $4 \leq f_j \leq 4s$ and $1 \leq j \leq s$, all the vertex labels will be distinct and from the set $\left[1,s\right] \cup \left[c+2,c+4s\right]$, except $0$, which will be repeated $s$ times. Considering the differences of these quadruples, we see that $\left(d_{j}+c\right) - j=c_{j}+c, \left(f_{j}+c\right)-j=e_{j}+c, d_{j}+c-0=d_{j}+c$, and $f_{j}+c-0=f_{j}+c$. Since, by construction, $c_{j},d_{j},e_{j},f_{j}$ are all distinct, we obtain $\left[c+1,c+4s\right]$ as the set of distinct edge labels.\end{pf}

Note that for the result in the first row of Table~\ref{2big} there are no restrictions on the sequence so the set of vertex labels is large. If we restrict the sequence, as in the second and third rows of Table~\ref{2big}, we will refine the set of vertex labels.
In Constructions~\ref{largea}--\ref{dddd}, we will give labellings of $C_{3}^{t}C_{4}^{s}$ when $2t < s \leq \frac{(13t+37)}{2}$, with the central vertex labelled $0$.

Let $P_{x}$ be the two-fold Skolem-type sequence of order $x$ given by Construction~\ref{like1}. Let $P_{x}^{'}$ be the the sequence obtained from $P_{x}$ by removing the pair $\left(1,1\right)$ from the end of the sequence. Define $P_{0}^{'}$ to be the sequence $\left(1,1\right)$.
%\textcolor[rgb]{1,0,0}{I think $s=i+l+j$}
%Let $L$ be a Langford sequence with defect $d$ and order $l=2d-1$.Let $S_t$ be a (hooked) Skolem sequence of order $t$.
%\textcolor[rgb]{1,0,0}{Note: I need to discuss that with Dr. Danny,
%\\ phase 0 for $s=2t+1$
%\\ phase 1 for $2t+2 \leq s \leq 3t+1$
%\\ phase 2x,2x+1 for $s \geq 3t+2$}
%\\
\begin{construction}\label{largea}
From a double Langford sequence with defect $d=t+1$ and order $l=2t+1$, construct the quadruples of the form $\left(0,d_{j}+c,j,f_{j}+c\right)$ as indicated in Construction~\ref{twofoldskolem03} with $c=t$. From a Skolem sequence of order $t$, construct the triples $\left(0,a_i+c,b_i+c\right)$ as indicated in Construction~\ref{skolem01} with $c=4l+t$. These triples and quadraples give a labelling for a $C_{3}^{t}C_{4}^{s}$ where $s=2t+1$, with the central vertex labelled $0$.
\end{construction}
%%%%%%%%%%%%%%%%%%%%%%%%%%%%%%%%%%%%%%%%%%%%%%%%%%%%%%%%%%%%%%%%%%%%%%%%%%%%%%%%%%%%%%%%%%%%%%%%%%%%%%%%%%%%%%%%
\begin{construction}\label{largeb}
By concatenating a double Langford sequence with defect $d=t+1$ and order $l=2t+1$, with a two-fold Skolem sequence of order $k$ ($k \leq t$), construct the quadruples of the form $\left(0,d_{j}+c,j,f_{j}+c\right)$ with $c=t$ as indicated in Construction~\ref{twofoldskolem03}. From a Skolem sequence of order $t$, construct the triples $\left(0,a_i+c,b_i+c\right)$ as indicated in Construction~\ref{skolem01} with $c=4k+4l+t$. These triples and quadruples give a labelling for a $C_{3}^{t}C_{4}^{s}$ where $2t+2 \leq s \leq 3t+1$, with the central vertex labelled $0$.
\end{construction}
%%%%%%%%%%%%%%%%%%%%%%%%%%%%%%%%%%%%%%%%%%%%%%%%%%%%%%%%%%%%%%%%%%%%%%%%%%%%%%%%%%%%%%%%%%%%%%%%%%%%%%%%%%%%%%%%
\begin{construction}\label{largec}
By concatenating a two-fold Skolem-type sequence $P_{x-1}^{'}$ of order $x-1$ and $1 \leq x \leq \frac{(t+3)}{2}$, with a double Langford sequence with defect $d=t+4x-1$ and order $l=2t+8x-3$, with the sequence $P_{0}^{'}$, and a two-fold Skolem-type sequence from Proposition~\ref{like3} of order $y$, for $0 \leq y \leq 4$, construct the quadruples of the form $\left(0,d_{j}+c,j,f_{j}+c\right)$ with $c=t$ as indicated in Construction~\ref{twofoldskolem03}. From a Skolem sequence of order $t$ construct the triples $\left(0,a_i+c,b_i+c\right)$ as indicated in Construction~\ref{skolem01} with $c=4l+t+4x+4y+4$. These triples and quadruples give a labelling for a $C_{3}^{t}C_{4}^{s}$ where $s=2t+9x+y-3$ and $3t+2 \leq s \leq \frac{(13t+37)}{2}$, with the central vertex labelled $0$.
\end{construction}
%%%%%%%%%%%%%%%%%%%%%%%%%%%%%%%%%%%%%%%%%%%%%%%%%%%%%%%%%%%%%%%%%%%%%%%%%%%%%%%%%%%%%%%%%%%%%%%%%%%%%%%%%%%%%%%%
\begin{construction}\label{dddd}
By concatenating a two-fold Skolem-type sequence $P_{x}$ of order $x$ with $1 \leq x \leq \frac{(t+3)}{2}$, with a double Langford sequence with defect $d=t+4x+1$ and order $l=2t+8x+1$, with a two-fold Skolem-type sequence from Proposition~\ref{like3} of order $y$, with $0 \leq y \leq 4$, construct the quadruples of the form $\left(0,d_{j}+c,j,f_{j}+c\right)$ with $c=t$ as indicated in Construction~\ref{twofoldskolem03}. From a Skolem sequence of order $t$ construct the triples $\left(0,a_i+c,b_i+c\right)$ as indicated in Construction~\ref{skolem01} with $c=4l+t+4x+4y$. These triples and quadruples give a labelling for a $C_{3}^{t}C_{4}^{s}$ where $s=2t+9x+y+1$ and $3t+2 \leq s \leq \frac{(13t+37)}{2}$, with the central vertex labelled $0$. \end{construction}%$s=9\left(\frac{x-1}{2}\right) +2t+1+y$
%%%%%%%%%%%%%%%%%%%%%%%%%%%%%%%%%%%%%%%%%%%%%%%%%%%%%%%%%%%%%%%%%%%%%%%%%%%%%%%%%%%%%%%%%%%%%%%%%%%%%%%%%%%%%%%%
Note that Constructions~\ref{largec} and~\ref{dddd} will cover all the values of $s \in \left[3t+2, \frac{(13t+37)}{2}\right]$.

The bound $s \leq \frac{(13t+37)}{2}$ comes from Constructions~\ref{largec} and~\ref{dddd} as we know $s=2t+9x+z$ or $x=\frac{(s-2t-z)}{9}$. Since $t \geq 2x-3$ then $s \leq \frac{(13t+2z+27)}{2}$ and since $z \in \left[-3,5\right]$ we obtain that $s \leq \frac{(13t+37)}{2}$.
\begin{lemma}\label{largee}
In the labellings of $C_{3}^{t}C_{4}^{s}$ given by Constructions~\ref{largea} -~\ref{dddd} for $4 \leq t \leq s \leq \frac{(13t+37)}{2}$, and $t \geq 2x-3$ when $s \geq 3t+2$,
\begin{enumerate}
	\item if $t \equiv 0,1\ ({\rm mod}\ 4)$, then the edge labels used are $\left[1,4s+3t\right]$ and the vertex labels used are from $\left\{0\right\} \cup \left[1,4s+3t\right]$, where each nonzero label occurs exactly once,
	\item if $t \equiv 2,3\ ({\rm mod}\ 4)$ then the edge labels used are $\left[1,4s+3t-1\right] \cup \left\{4s+3t+1\right\}$ and the vertex labels used are from $\left\{0\right\} \cup \left[1,4s+3t-1\right]  \cup \left\{4s+3t+1\right\}$, where each nonzero label occurs exactly once.	
\end{enumerate}
\end{lemma}
%Using a two-fold Skolem-like sequence $C^{1}$ with Construction~\ref{twofoldskolem03}
\begin{pf}We will prove these results for Construction~\ref{dddd}. The proof for Construction~\ref{largec} follows in the same fashion as Construction~\ref{dddd}. For Constructions~\ref{largea} and~\ref{largeb} the proofs also follow in the same fashion as Construction~\ref{dddd}, but with no Skolem-type sequences.

Let $t \equiv 0,1\ ({\rm mod}\ 4)$. We begin by considering the edge labels.
%%%%%%%%%%%%%%%%%%%%%%%%%%%%%%%%%%%%%%%%%%%%%%%%%%%%%%%%%%%%%%%%%%%%%%%%%%%%%%%%%%%%%%%%%%%%%%
Consider the quadruples formed by the concatenated sequence in Construction~\ref{dddd}, with $c=t$. Those quadruples corresponding to $P_{x}$ yield edge labels $\left[t+1,t+4x\right]$, by Table~\ref{2big} (row 7). For the quadruples corresponding to the double Langford sequence, we obtain the edge labels $\left[t+4x+1,4l+t+4x\right]$, by Table~\ref{2big} (row 4), by considering $c=t+4x$ (the length of $P_{x}$). For the quadruples corresponding to the two-fold Skolem-type sequence from Proposition~\ref{like3}, we obtain the edge labels $\left[4l+t+4x+1,4l+t+4x\right. \\
    \left.+4y\right]$, by Table~\ref{2big} (rows 9-12), by considering $c=4l+t+4x$ (the length of $P_{x}$ and the double Langford sequence). Note that if $y=0$, then we are considering $C^{0}$, the empty sequence, and hence produce no edge labels.
%%%%%%%%%%%%%%%%%%%%%%%%%%%%%%%%%%%%%%%%%%%%%%%%%%%%%%%%%%%%%%%%%%%%%%%%%%%%%%%%%%%%%%%%
% For a two-fold Skolem-like sequence $P_{x}$ of order $x$ with $c=0$, by Table~\ref{2big} (row 7) this construction yields edge labels $\left[1+t,4x+t\right]$.\\
%For a double Langford sequence defect with defect $d=t+4x+1$ and order $l=2t+8x+1$ and $c=4x$ (we concatenated this sequence to a two-fold Skolem-like sequence $P_{x}$ of order $x$ then the shift $c=4x$), by Table~\ref{2big} (row 4) this construction yields edge labels $\left[1+t+4x,4l+4x+t\right]$.\\
%For a two-fold Skolem-like sequence from Proposition~\ref{like3} of order $y$ with $c=4x+4l$ (we concatenated this sequence to a two-fold Skolem-like sequence $P_{x}$ of order $x$ and a double Langford sequence of order $l$ then the shift $c=4x+4l$), by Table~\ref{2big} (row 9-12) this construction yields edge labels\\ $\left[1+t+4x+4l,4y+t+4l+4x\right]$. Note that if $y=0$, then we are considering $C^{0}$, the empty set, and hence produce no edge labels.\\

Consider the triples formed by the Skolem sequence of order $t$ with $c=4l+t+4x+4y$ (the length of $P_{x}$, the double Langford sequence, and the two-fold Skolem-type sequence). By Lemma~\ref{l1}(1), this construction yields edge labels $\left[1,t\right] \cup \left[4l+t+4x+4y+1,4l+3t+4x+4y\right]$.

From the above discussion, it is clear that all edges are distinct and are exactly the set $\left[1,4s+3t\right]$, where $s=l+x+y$.
%%%%%%%%%%%%%%%%%%%%%%%%%%%%%%%%%%%%%%%%%%%%%%%%%%%%  vertecise  %%%%%%%%%%%%%%%%%%%%%%%%%%%%%%%%%%%%%%%%%%%%%%%%%%%%%%%%%%%%%%%%%%%%%%%%%%%%%%%%%%%%%%%%%%%%%%%%%%%%%%%

%%%%%%%%%%%%%%%%%%%%%%%%%%%%%%%%%%%%%%%%%%%%%%%%%%%%%%%%%%%%%%%%%%%%%%%%%%%%%%%%%%%%%%%%%%%%%%
We now consider the vertex labels. Consider the quadruples formed by the concatenated sequence in Construction~\ref{dddd}, with $c=t$. Those quadruples corresponding to $P_{x}$ yield vertex labels that are a subset of $\left\{0,1\right\} \cup \left[t+2x,t+4x\right] \cup \left\{4i | 1 \leq i \leq x-1\right\}$, by Table~\ref{2big} (row 7) and since $t \geq 2x-3$ there are no vertices repeated.
For the quadruples corresponding to the double Langford sequence, we obtain vertex labels from $\left\{0\right\} \cup \left[t+4x+1,9t+36x+4\right]$, by Table~\ref{2big} (row 4), by considering $c=t+4x$ (the length of $P_{x}$). For the quadruples corresponding to the two-fold Skolem-type sequence from Proposition~\ref{like3}, we obtain vertex labels from $\left\{0,2,3,5,6\right\} \cup \left[9t+36x+7,9t+36x+4y+4\right]$, by Table~\ref{2big} (rows 9-12), by considering $c=4l+t+4x$ (the length of $P_{x}$ and double Langford sequence). The labels $\left\{0,2,3,5,6\right\}$ are only used once so they do not conflict with any other vertex labels, since $t+2x \geq 6$ and none of these labels is a multiple of four. Note that if $y=0$, then we are considering $C^{0}$, the empty sequence, and hence produce no vertex labels.
%%%%%%%%%%%%%%%%%%%%%%%%%%%%%%%%%%%%%%%%%%%%%%%%%%%%%%%%%%%%%%%%%%%%%%%%%%%%%%%%%%%%%%%%
%For a two-fold Skolem-like sequence of order $x$ with $c=0$, by Table~\ref{2big} (row 7) this construction yields vertex labels from is a set of $x$ distinct positive integers.\\
%From Table~\ref{tablek}, we can construct a two-fold Skolem-like sequence construction of order $n$, where $H=\left\{\left\{1\right\} \cup \left\{4i\right\} | 1 \leq i \leq n-1\right\}$ is a set of $n$ distinct positive integers.
%For a double Langford sequence defect with defect $d=t+4x+1$ and order $l=2t+8x+1$ and $c=4x$, by Table~\ref{2big} (row 4) this construction yields vertex labels from $\left\{0\right\} \cup \left[t+4x+1,3t+12x+1\right] \cup \left[3t+12x+2,9t+36x+4\right]$.\newline
%For a two-fold Skolem-like sequence from Proposition~\ref{like3} of order $y$ with $c=4x+4l$, by Table~\ref{2big} (row 9-12) this construction yields vertex labels $\left\{0,2,3,5,6\right\} \cup \left[9t+36x+7,9t+36x+4y+4\right]$. Note that if $y=0$, then we are considering $C^{0}$, the empty set, and hence produce no vertex labels.

Consider the triples formed by the Skolem sequence of order $t$ with $c=4l+t+4x+4y$ (the length of $P_{x}$, double Langford sequence, and the two-fold Skolem-type sequence). By Lemma~\ref{l1}(1), this construction yields vertex labels from $\left\{0\right\} \cup \left[9t+36x+4y+5,11t+36x+4y+4\right]$.
%For a Skolem sequence of order $t$ construct the triples with $c=4x+4l+4y$, then by Lemma~\ref{l1} (1), this construction yield vertex labels from $\left\{0\right\} \cup \left[9t+36x+4y+5,11t+36x+4y+4\right]$.\\

From the above discussion, it is clear that all vertex labels are distinct and are from the set $\left[0,4s+3t\right]$ where $s=l+x+y$.

If $t \equiv 2,3\ ({\rm mod}\ 4)$, we proceed similarly to the proof for the case $t \equiv 0,1\ ({\rm mod}\ 4)$, but use a hooked Skolem sequence instead of Skolem sequence with Construction~\ref{skolem01} and Lemma~\ref{l1}(2).\end{pf}
\section{$C_{3}^{t}C_{4}^{s}$}\label{section4}
In this section, we prove (near) graceful labellings exist for $C_{3}^{t}C_{4}^{s}$.
%\\we have to check the geacefulness \\
%no labels show up twice\\
%vertex label from ${0,\ldots,m}$, show up once\\
%edge label from ${1,\ldots,m}$, show up once\\
%\subsection{$t \geq s$}
\begin{theorem}\label{firstthm}If $G=C_{3}^{t}C_{4}^{s}$, where $t \geq s \geq 1$, then $G$ is graceful when $t \equiv 0,1\ ({\rm mod}\ 4)$ and near graceful when $t \equiv 2,3\ ({\rm mod}\ 4)$.\end{theorem}
\begin{pf}
\textbf{Case 1}: $t \equiv 0,1\ ({\rm mod}\ 4)$.

%Start by a two-fold Skolem sequence or order $s$ and it is exists by Table~\ref{sbig}.
Use a two-fold Skolem sequence of order $s$, with Construction~\ref{twofoldskolem03} and with $c=t$, to get $\left(0,d_{j}+c,j,f_{j}+c\right)$, $1\leq j\leq s$. Using a Skolem sequence of order $t$ in Construction~\ref{skolem01} with $c=4s+t$ and $t \geq s$, gives $\left(0,a_i+c,b_i+c\right)$, $1\leq i\leq t$. These sequences are known to exist; see Table~\ref{sbig}.

These two constructions give the vertex labels and the induced edge labels of $G$. By Table~\ref{2big} (row 1) the quadruples use vertex labels in $\left[0,s\right] \cup \left[t+2,4s+t\right]$ and edge labels $\left[t+1,4s+t\right]$ and by Lemma~\ref{l1}(1) the triples use vertex labels in $\left\{0\right\} \cup \left[4s+t+1,4s+3t\right]$ and edge labels $\left[1,t\right] \cup \linebreak \left[4s+t+1,4s+3t\right]$. Thus, this is a graceful labelling.

\textbf{Case 2}: $t \equiv 2,3\ ({\rm mod}\ 4)$.\\
This is similar to Case 1 using a hooked Skolem sequence and Lemma~\ref{l1}(2) instead of a Skolem sequence and Lemma~\ref{l1}(1).\end{pf}
As an example, consider the two-fold Skolem sequence $\left(3,1,1,3,2,2,2,2\right.$, $\left.3,1,1,3\right).$ This sequence gives quadruples and triples that together gracefully label $G=C_{3}^{4}C_{4}^{3}$ (see Figure~\ref{ga}).
\begin{theorem}\label{secondthm}If $G=C_{3}^{t}C_{4}^{s}$, where $4 \leq t \leq s \leq \frac{(13t+37)}{2}$, then $G$ is graceful when $t \equiv 0,1\ ({\rm mod}\ 4)$ and is near graceful when $t \equiv 2,3\ ({\rm mod}\ 4)$.\end{theorem}
\begin{pf}
\noindent\textbf{Case 1}: $t \equiv 0,1\ ({\rm mod}\ 4)$. For $4 \leq t < s \leq 2t $ and $s$ odd, use a two-fold Skolem sequence of order $s$ as given in Table~\ref{table2fo}, with Construction~\ref{twofoldskolem03} and with $c=t$, to get $\left(0,d_{j}+c,j,f_{j}+c\right)$, $1\leq j\leq s$. Using a Skolem sequence of order $t$ in Construction~\ref{skolem01} with $c=4s+t$ and $t < s$, gives $\left(0,a_i+c,b_i+c\right)$, $1\leq i\leq t$. These sequences exist, as detailed in Table~\ref{sbig}.

These two constructions give the vertex labels and the induced edge labels of $G$. By Table~\ref{2big} (row 2) the quadruples use vertex labels in $\left[0,s\right] \cup \left[\frac{(s+3)}{2}+t,4s+t\right]$ and edge labels $\left[t+1,4s+t\right]$ and by Lemma~\ref{l1}(1) the triples use vertex labels in $\left\{0\right\} \cup \left[4s+t+1,4s+3t\right]$ and produce edge labels $\left[1,t\right] \cup \left[4s+t+1,4s+3t\right]$. Thus, this is a graceful labelling. %done
%%%%%%%%%%%%%%%%%%%%%%%%%%%%%%%%%%%%%%%%%%%%%%%%%%%%

For $4 \leq t < s \leq 2t $ and $s$ even, use a two-fold Skolem sequence of order $s$ as given in Table~\ref{table2fe}, with Construction~\ref{twofoldskolem03} and $c=t$, to get $\left(0,d_{j}\!+\!c,j,f_{j}\!+\!c\right)$, $1\leq j\leq s$. Using a Skolem sequence of order $t$ in Construction~\ref{skolem01} with $c=4s+t$ and $t < s$, gives $\left(0,a_i+c,b_i+c\right)$, $1\leq i\leq t$. These sequences exist, as detailed in Table~\ref{sbig}.

These two constructions give the vertex labels and the induced edge labels of $G$. By Table~\ref{2big} (row 2) the quadruples use vertex labels in $\left[0,s\right] \cup \left[\frac{(s+2)}{2}+t,4s+t\right]$ and edge labels $\left[t+1,4s+t\right]$ and by Lemma~\ref{l1}(1) the triples use vertex labels in $\left\{0\right\} \cup \left[4s+t+1,4s+3t\right]$ and edge labels $\left[1,t\right] \cup \left[4s+t+1,4s+3t\right]$. Thus, this is a graceful labelling.

For $4 \leq t < s=2t+1 $ use Construction~\ref{largea}, and for $2t+2 \leq s \leq 3t+1$ use Construction~\ref{largeb}. These constructions give the vertex labels and the induced edge labels of $G$. By Lemma~\ref{largee}(1) the quadruples and triples use vertex labels in $\left[0,4s+3t\right]$ and edge labels $\left[1,4s+3t\right]$. Thus, this is a graceful labelling.

We now consider the case $3t+2 \leq s \leq \frac{(13t+37)}{2}$. Define $I_{x}=[2t+9x-3$, $2t+9x+1]$, and $J_{x}=\left[2t+9x+1,2t+9x+5\right]$ for $x \geq 1$ and fixed $t$. Define $K_{x}=\left[2t+9x-3, 2t+9x+5\right]=I_{x} \cup J_{x}$. Note that $K_{x} \cap K_{x+1} = \emptyset$, but for any interval $K_{x}$ the largest element is $2t+9x+5$ and the smallest element in $K_{x+1}$ is $2t+9x+6$. Note $\bigcup_{x \geq 1} K_{x}=\left[2t+6,\infty\right)$, and so for all integers $s \in \left[2t+6,\infty\right)$, there exists $x$ such that $s \in K_{x}$, and hence $s \in I_{x}$ or $s \in J_{x}$. That is, $s$ can be written either in the form $s=2t+9x+y-3$ or the form $s=2t+9x+y+1$, where $0 \leq y \leq 4$, for some $x$.% (since $t \geq 4$ and $x \geq 1$ therefore $2t+9x-3 \geq 3t+2$)

%%%%%%%%%%%%%%%%%%%%%%%%%%%%%%%%%%%%%%%%%%%%%%%%%% large c
For $s=2t+9x+y-3$, and $3t+2 \leq s \leq \frac{(13t+37)}{2}$ use Construction~\ref{largec}, and for $s=2t+9x+y+1$, and $3t+2 \leq s \leq \frac{(13t+37)}{2}$ use Construction~\ref{dddd}. These constructions give the vertex labels and the induced edge labels of $G$. By Lemma~\ref{largee}(1) the quadruples and triples use vertex labels in $\left[0,4s+3t\right]$ and edge labels $\left[1,4s+3t\right]$. Thus, this is a graceful labelling.

%%%%%%%%%%%%%%%%%%%%%%%%%%%%%%%%%%%%%%%%%%%%%%%%%% large d
%For $s=9x+2t+1+y,$ use Construction~\ref{dddd}. This construction give the vertex labels and the induced edge labels of $G$. By Lemma~\ref{largee}(1) the quadruples and triples use vertex labels in $\left\{0\right\} \cup \left[1,3t+4s\right]$ and edge labels $\left[1,3t+4s\right]$. Thus, this is a graceful labelling.
%These two constructions give the labels of $G$. By Table~\ref{2big} (row 4) the quadruples use vertex labels in $\left\{0\right\} \cup \left\{t+1,\ldots,4s+t\right\}$ and edge labels $\left\{t+1,\ldots,t+4s\right\}$ and by Lemma~\ref{l1}(1) the triples use vertex labels in $\left\{0\right\} \cup \left\{1+t+4s,\ldots,3t+4s\right\}$ and edge labels $\left\{1,\ldots,t\right\} \cup \left\{1+t+4s,\ldots,3t+4s\right\}$. Thus, this is a graceful labelling.\\
%\textcolor[rgb]{1,0,0}{I will do the rest after the confirmation ?} \\
%\indent For $2t < s $, use the appropriate construction from Constructions~\ref{large0} - \ref{large2x+1}. These two constructions give the labels of $G$. By Lemma~\ref{largee} (1) we use edge labels $\left\{1,2,\ldots,3t+4s\right\}$ and vertex labels $\left\{0\right\} \cup \left\{1,\ldots,3t+4s\right\}$, when $t \equiv 0\ or\ 1\ ({\rm mod}\ 4)$.\\
%Finally, we can get gracefully labelled graph $G$ with $s > t \geq 4$ which uses edge labels $\left\{1,2,\ldots,3t+4s\right\}$ and vertex labels from $\left\{0,1,\ldots,3t+4s\right\}$, and no repeated edge or vertex labels. So $G=C_{3}^{t}C_{4}^{s}$, where $s > t \geq 4$, is graceful when $t \equiv 0,1\ ({\rm mod}\ 4)$.\\

\noindent\textbf{Case 2}: $t \equiv 2,3\ ({\rm mod}\ 4)$. Similar to Case 1, but use a hooked Skolem sequence, Lemma~\ref{l1}(2), and Lemma~\ref{largee}(2) instead of a Skolem sequence, Lemma~\ref{l1}(1) and Lemma~\ref{largee}(1). Then, this is a near graceful labelling.\end{pf}

In Theorem~\ref{firstthm} and Theorem~\ref{secondthm} we proved (near) graceful labellings exist for $C_{3}^{t}C_{4}^{s}$ if $1 \leq s \leq t$ or $4 \leq t \leq s \leq \frac{(13t+37)}{2}$, omitting the cases for $t=1,2,3$. So, in the following lemma we consider those cases.

\begin{lemma}\label{BZ}For $k \geq 1$ and $x \ge 1$, if
\begin{enumerate}
\item $x \equiv 0 \!\!\pmod{4}$ when $w \geq 1$ and $\frac{(2k-12w+2)}{4} \leq s \leq \frac{(6k-12w-5)}{4}$, or
\item $x \equiv 1 \!\!\pmod{4}$ when $w \geq 0$ and $\frac{(2k-12w-1)}{4} \leq s \leq \frac{(6k-12w-8)}{4}$, and
\end{enumerate}
a graceful labelling of $C_3^{4w+x}C_4^{s}$ exists, then a graceful labelling of \\$C_3^{4w+x}C_4^{s+4k-1}$ exists. Alternatively, if
\begin{enumerate}
\item[3.] $x \equiv 2 \!\!\pmod{4}$ when $w \geq 0$ and $\frac{(2k-12w-2)}{4} \leq s \leq \frac{(6k-12w-9)}{4}$, or
\item[4.] $x \equiv 3 \!\!\pmod{4}$ when $w \geq 0$ and $\frac{(2k-12w-3)}{4} \leq s \leq \frac{(6k-12w-10)}{4}$, and
\end{enumerate}
a near graceful labelling of $C_3^{4w+x}C_4^{s}$ exists which contains a triangle $(0$, $4s+12w+5$, $4s+12+7)$ if $x \equiv 2 \!\!\pmod{4}$ or contains triangles $(0$, $4s+12w+7$, $4s+12w+8)$ and $(0$, $4s+12w+6$, $4s+12w+10)$, then a near graceful labelling of $C_3^{4w+x}C_4^{s+4k-1}$ exists.
\end{lemma}

\begin{pf} First, note that if $C_3^{4w+x}C_4^{s}$ is graceful, it uses edge labels $[1$, $4s+12w+3x]$ and vertex labels from $\left[0,4s+12w+3x\right]$. If if $C_3^{4w+x}C_4^{s}$ is near graceful, it uses edge labels $\left[1,4s+12w+3x-1\right] \cup \{4s+12w+3x+1\}$ and vertex labels from $\left[0,4s+12w+3x-1\right] \cup \{4s+12w+3x+1\}$. We consider four cases corresponding to the four possibilities of the lemma.

\noindent{\bf Case 1:} By using a two-fold Langford sequence with defect $d=6k-1$ and order $l=4k-1$ from Table~\ref{table2l} with Construction~\ref{twofoldskolem03}, $c=4s+12w$, we can label a $C_4^{4k-1}$ with edge labels $\left[4s+12w+1,16k+4s+12w-4\right]$ and vertex labels from $\left\{0\right\} \cup \left[6k-1,10k-3\right] \cup \left[8k+4s+12w-1,16k+4s+12w-4\right]$.

Note that to make sure any vertex label is only used at most once, $6k-1$ must be greater than $4s+12w$ and since $\frac{(6k-12w-5)}{4} \geq s$, we have no vertex labels used twice. Also, to avoid a conflict in vertex labelling, $8k+4s+12w-1 > 10k-3$ and since $s \geq \frac{(2k-12w+2)}{4}$, we have no vertex labels used twice.

By identifying the vertices with label zero in $C_3^{4w}C_4^{s}$ and $C_4^{4k-1}$, the labelling we obtain is a graceful labelling of $C_3^{4w}C_4^{s+4k-1}$ with edge labels $\left[1,16k+4s+12w-4\right]$ and vertex labels from $\left[0,16k+4s+12w-4\right]$.

\noindent{\bf Case 2:} Again, by using a two-fold Langford sequence with defect $d=6k-1$ and order $l=4k-1$ from Table~\ref{table2l} with Construction~\ref{twofoldskolem03}, $c=4s+12w+3$, we can label a $C_4^{4k-1}$ with edge labels $\left[4s+12w+4,16k+4s+12w-1\right]$ and vertex labels from $\left\{0\right\} \cup \left[6k-1,10k-3\right] \cup [8k+4s+12w+2$, $16k+4s+12w-1]$. As in the previous case, by choice of $s$, we avoid re-using labels.

Then by identifying the vertices with label zero in $C_3^{4w+1}C_4^{s}$ and $C_4^{4k-1}$, we obtain is a graceful labelling of $C_3^{4w+1}C_4^{s+4k-1}$.

\noindent{\bf Case 3:}
By using a two-fold Langford sequence with defect $d=6k-1$ and order $l=4k-1$ from Table~\ref{table2l} with Construction~\ref{twofoldskolem03}, $c=4s+12w+4$, we can label a $C_4^{4k-1}$ with edge labels $\left[4s+12w+5,16k+12w+4s\right]$ and vertex labels from $\left\{0\right\} \cup \left[6k-1,10k-3\right] \cup \left[8k+4s+12w+3,16k+4s+12w\right]$. As in case 1, choice of $s$ allows us to avoid re-using vertex and edge labels.

If we replace the triangle containing edge length $2$, $(0$, $4s+12w+5$, $4s+12w+7)$, by $\left(0\right.$, $16k+4s+12w+1$, $16k+4s+12w+3)$ and identify the vertices with label zero in the $C_3^{4w+2}C_4^{s}$ and the $C_4^{4k-1}$ labelling, then we obtain a near graceful labelled $C_3^{4w+2}C_4^{s+4k-1}$.

\noindent{\bf Case 4:} By using a two-fold Langford sequence with defect $d=6k-1$ and order $l=4k-1$ from Table~\ref{table2l} with Construction~\ref{twofoldskolem03}, $c=4s+12w+5$, we can label a $C_4^{4k-1}$ with edge labels $\left[4s+12w+6,16k+4s+12w+1\right]$ and vertex labels from $\left\{0\right\} \cup \left[6k-1,10k-3\right] \cup \left[8k+4s+12w+4,16k+4s+12w+1\right]$. As in case 1, by choice of $s$ we eliminate the possibility of re-using labels.

Replace the triangles containing edge lengths $1$ and $4$, $(0$, $4s+12w+7$, $4s+12w+8)$ and $(0$, $4s+12w+6$, $4s+12w+10)$, by $(0$, $16k+4s+12w+3$, $16k+4s+12w+4)$ and $(0$, $16k+4s+12w+2$, $16k+4s+12w+6)$. If we then identify the vertices with label zero in the $C_3^{4w+3}C_4^{s}$ and the $C_4^{4k-1}$ labelling, we obtain a near graceful labelled $C_3^{4w+3}C_4^{s+4k-1}$.
\end{pf}

Note that in the near graceful cases, the resulting near graceful labellings contain triangles of the appropriate form to satisfy the criteria of the lemma, allowing an iterative use.

One of the difficulties of Lemma~\ref{BZ} is that it gives what is essentially an existential result. Consider the graph $C_3^4C_4^{100}$. Then certainly, if we are going to apply Lemma~\ref{BZ}, we must consider its first case. But for what combination of $s$ and $k$? We must find a pair $(s,k)$ such that $s+4k-1 = 100$ and such that $\frac{(2k-10)}{4} \leq s \leq \frac{(6k-117)}{4}$. One such pair is $(21,20)$. In that case, if $C_3^4C_4^{21}$ is graceful, we are done. Happily, we can conclude that $C_3^4C_4^{21}$ is graceful, by Theorem~\ref{secondthm}. Of course, this isn't the only pair that satisfies these two constraints, and different pairs will produce different labellings. (The pair $(17,21)$ is another such solution, since Theorem~\ref{secondthm} tells us that $C_3^4C_4^{17}$ is graceful.)

Thus, to prove that an arbitrary graph $C_3^{t}C_4^{s}$ is (near) graceful, we will show that an appropriate $k$ can be found to satisfy the appropriate condition in Lemma~\ref{BZ}, and that a smaller (near) graceful labelling exists by one of Theorems~\ref{firstthm} or \ref{secondthm}.

\begin{theorem}\label{C1}
If $s\in \mathbb{N}$, and $t=4w+z$ with $z=0,1,2,3$ then $C_3^{t}C_4^{s}$ can be gracefully labelled if $t=4w,4w+1$ and near gracefully labelled if $t=4w+2,4w+3$.\end{theorem}

\begin{pf} We consider four cases based on the modularity of $t$. We will repeatedly use induction on $s$, the number of $4$-cycles in $C_3^{4w+z}C_4^{s}$. For $w=0$ and $z = 1$, 2, or 3, see Appendix B for (near) graceful labellings when $s$ is small. For $w \ge 1$ and any $z$ note that the base cases needed for the induction (below the indicated value of $s$ in each case) can be given by Theorems~\ref{firstthm} and \ref{secondthm}.

\noindent\textbf{Case 1}: $t=4w$. Define $I_{k}$ to be the real interval $\left[\frac{(18k-12w-2)}{4},\frac{(22k-12w-9)}{4}\right]$ for fixed $k$. If $k \geq \frac{25}{4}$ then $\frac{(18\left(k+1\right)-12w-2)}{4} \leq \frac{(22k-12w-9)}{4}$. That is, $I_{k} \cap I_{k+1} \neq \emptyset$. This implies $\bigcup_{k \geq 7} I_{k}=\left[\frac{124-12w}{4},\infty\right)$, for fixed arbitrary $w$, and so for all $s \in \left[\frac{121-12w}{4},\infty\right),$ there exists $k$ such that $s \in I_{k}.$

We proceed by induction on $s$. Let $s \geq 28$ be an integer. There exists some $k$ such that $s \in I_{k}$. Therefore, letting $s=4k+s'-1$, then $\frac{(2k-12w+2)}{4} \leq s' \leq \frac{(6k-12w-5)}{4}$. By induction, a graceful labelling of $C_3^{4w}C_4^{s'}$ exists, and by Lemma~\ref{BZ} a graceful labelling of $C_3^{4w}C_4^{s}$ exists.

\noindent\textbf{Case 2}: $t=4w+1$. As in case 1, if $I_{k}= \left[\frac{(18k-12w-5)}{4},\frac{(22k-12w-12)}{4}\right]$, then for all $s \in \bigcup_{k \geq 7} I_{k} = \left[\frac{121-12w}{4},\infty\right),$ there exists $k$ such that $s \in I_{k}.$ Again, following by induction on $s$, a graceful labelling of $C_3^{4w+1}C_4^{s'}$ exists, and by Lemma~\ref{BZ} a graceful labelling of $C_3^{4w+1}C_4^{s}$ exists.

\noindent\textbf{Case 3}: $t=4w+2$.  Again, define $I_{k} = \left[\frac{(18k-12w-6)}{4},\frac{(22k-12w-13)}{4}\right]$. Then for all $s \in \bigcup_{k \geq 7} I_{k} = \left[\frac{120-12w}{4},\infty\right),$ there exists $k$ such that $s \in I_{k}.$ With induction on $s$, a graceful labelling of $C_3^{4w+2}C_4^{s'}$ exists, and by Lemma~\ref{BZ} a graceful labelling of $C_3^{4w+2}C_4^{s}$ exists.

\noindent\textbf{Case 4}: $t=4w+3$. As previously, let $I_{k} = \left[\frac{(18k-12w-7)}{4},\frac{(22k-12w-14)}{4}\right]$, so that for all $s \in \bigcup_{k \geq 7} I_{k} = \left[\frac{119-12w}{4},\infty\right),$ there exists $k$ such that $s \in I_{k}.$ Again using induction on $s$, a graceful labelling of $C_3^{4w+3}C_4^{s'}$ exists, and by Lemma~\ref{BZ} a graceful labelling of $C_3^{4w+3}C_4^{s}$ exists.
\end{pf}

\section{$C_{3}^{t}C_{5}^{p}$}\label{section5}
Using Theorem~\ref{thirdthm} with Langford sequences, we can obtain (near) graceful labellings of $C_{3}^{t}C_{5}^{p}$.

In Constructions~\ref{c35}, we will give labellings of $C_{3}^{t}C_{5}^{p}$, with the central vertex labelled $0$.
\begin{construction}\label{c35}
Given a $C_{5}^{p}$ gracefully labelled by the 5-tuples $(0$, $d_{b_{i}}+p$, $b_i$, $a_i$, $d_{a_{i}}+p)$ for each $1 \leq i \leq p$ formed by
\begin{enumerate}
\item \label{c35.0}Construction~\ref{c5a}(1), if $p \equiv 0 \!\!\pmod{4}$, or
\item \label{c35.1}Construction~\ref{c5a}(2), if $p \equiv 1 \!\!\pmod{4}$, or
\item \label{c35.2}Construction~\ref{c5a}(3), if $p \equiv 2 \!\!\pmod{4}$, or
\item \label{c35.3}Construction~\ref{c5a}(4), if $p \equiv 3 \!\!\pmod{4}$, then
\end{enumerate}
replace the 5-tuples with $\left(0,d_{b_{i}}+p+3t,b_i,a_i,d_{a_{i}}+p+3t\right)$ for each $1 \leq i \leq p$ and $t \geq 2p+1$. From a Langford sequence with defect $p+1$ and order $t$, form triples via Construction~\ref{skolem01} with $c=p+t$. These 5-tuples and triples give a labelling of a $C_{3}^{t}C_{5}^{p}$, with the central vertex labelled $0$.
%give a labelling of a Ct3 windmill, with the central vertex labelled 0.
\end{construction}
Recall the result of Theorem~\ref{thirdthm}: if $G=C_{5}^{p}$, then $G$ is near graceful when $p \equiv 1,2\ ({\rm mod}\ 4)$ and $G$ is graceful when $p \equiv 0,3\ ({\rm mod}\ 4)$. We consider both of these cases in the following theorem.%In order to apply Construction~\ref{c35} we have two cases: $p \equiv 0\ ({\rm mod}\ 4)$ and $p \equiv 3\ ({\rm mod}\ 4)$.% and Construction~\ref{c35n} we have two cases $p \equiv 1\ ({\rm mod}\ 4)$ and $p \equiv 2\ ({\rm mod}\ 4)$.

\begin{lemma}\label{c35c}The labellings of $C_{3}^{t}C_{5}^{p}$ given by Construction~\ref{c35}(\ref{c35.0}) and \ref{c35}(\ref{c35.3})  use edge labels $\left[1,5p+3t\right]$ exactly once and vertex labels from $\left[0,5p+3t\right]$, each nonzero label occurring exactly once, whereas the labellings given by Construction~\ref{c35}(\ref{c35.1}) and \ref{c35}(\ref{c35.2}) use edge labels $\left[1,5p+3t-1\right] \cup \left\{5p+3t+1\right\}$ exactly once and vertex labels from $\left[0,5p+3t-1\right]  \cup \left\{5p+3t+1\right\}$, each nonzero label occurring exactly once.
\end{lemma}
\begin{pf}We have four cases based on congruence modulo 4. We begin with $p \equiv 0\ ({\rm mod}\ 4)$ and $p \equiv 3\ ({\rm mod}\ 4)$ where we are using gracefully labelled $C_{5}^{p}$. %We will prove one case only. The other cases follow similarly.\\

\noindent\textbf{Case 1:} If $p \equiv 0\ ({\rm mod}\ 4)$, form the $5$-tuples $(0$, $d_{b_{i}}+p+3t$, $b_i$, $a_i$, $d_{a_{i}}+p+3t)$ for each $1 \leq i \leq p$. We begin by considering vertex labels.
		
% Start with a Skolem sequence of order $p=4m$ and construct the pairs $\left(a_i,b_i\right)$ with $1 \leq i \leq p$ as indicated in Construction~\ref{c5a}. Also, form a Skolem sequence of order $2p$ from Table~\ref{table1} and construct the pairs $\left(c_j,d_j\right)$ with $1 \leq j \leq 2p$ as indicated in Construction~\ref{c5a}.%For the $d_{a_i},d_{b_i} \geq p+1$ so the minimum value of the second and fifth entries of the $5$-tuples is $2p+3t+1$ then all the entries of the $5$-tuples are distinct.\\
In a Skolem-type sequence, the $a_i,b_i,c_i$ and $d_i$ are all unique. From the Skolem sequence $S_{1}$ of order $p$ indicated in Construction~\ref{c5a}, we notice that the entries $a_i$ and $b_i$ in the third and fourth entries of the $5$-tuples give the distinct numbers $\left[1,2p\right]$. Recall that the Skolem sequence $S_{2}$ of order $2p$ used in Construction~\ref{c5a} contains no right endpoints in the first $p$ positions, with the first right endpoint occurring in position $p+1$. Note $\bigcup\limits_{i=1}^{p} \left\{d_{a_i},d_{b_i}\right\}$ is the set of all right endpoints, and so a subset of $\left[p+1,4p\right].$ From the $5$-tuples, we use the right endpoints as indicated in the second and fifth entries. We add $p+3t$ to each element based on our construction so we obtain a subset of $\left[2p+3t+1,5p+3t\right]$.
	
Consider the triples of the form $\left(0,a_{i}+p+t,b_{i}+p+t\right)$ with $1 \leq i \leq t$ given by using Construction~\ref{skolem01} with $c=p+t$. These triples give the vertex labels for $t$ triangles $C_{3}^{t}$, where $0$ is repeated $t$ times as a common vertex. Thus by Lemma~\ref{llang1} these triples use vertex labels from $\left\{0\right\} \cup \left[p+t+1,p+3t\right]$.

Thus, in this labelling, the common vertex $0$ is repeated $p+t$ times. The remaining vertices are all distinct and from the union of disjoint sets $\left[1,2p\right] \cup \left[p+t+1,p+3t\right] \cup \left[2p+3t+1,5p+3t\right]$.

%%%%%%%%%%%%%%%%%%%%%%%%%%%%%%%%%%%%%%%%%%%%%%%%%%%%%%%%%%%%%%%%%%%%%%%%%%%%%%%%%%%%%%%%%%%%%%%%%%%%%%%%%%%%%%%%%%%%%%%%%%%%%%%%%%%%%%%%%%%%%%%%%%%%%%%%%%%%%%%%
%%%%%%%%%%%%%%%%%%%%%%%cccccccccccccccccccccccccccccccccccccccccccccccccccccccccccccccccccccccccccccccccccccccccccccccccccccccccccccccccccccccccccccccccccc
We now examine the edge labels from this construction, by considering the differences between subsequent entries (taken cyclically) in $\left(0\right.$, $d_{b_i}+p+3t$, $b_i$, $a_i$, $d_{a_i}+p+3t)$. The differences between $b_i$ and $a_i$ produce the distinct numbers $\left[1,p\right].$
		
The differences $\left(d_{b_i}+p+3t\right)$ $-$ $0$ and the differences $\left(d_{a_i}+p+3t\right) - 0$ produce distinct numbers. Call this set of distinct numbers $A$ and observe that $A \subseteq \left[2p+3t+1,5p+3t\right]$. Also, from $\left(d_{b_i} + p+3t\right) - b_i= \left(d_{b_i} - b_i\right) +p+3t = c_{b_i} + p+3t$ and $\left(d_{a_i} + p+3t\right) - a_i= \left(d_{a_i} - a_i\right) +p+3t = c_{a_i} + p+3t$, we get distinct numbers. Call this set of distinct numbers $B$ and observe that $B \subseteq \left[p+3t+1,5p+3t-1\right]$. We know in a Skolem sequence, the $c_i$ and $d_i$ are all unique so $A$ and $B$ are disjoint. Note that $\left|A \cup B\right|=4p$. Now, we can conclude that all the previous differences give exactly the edge labels $\left[p+3t+1,5p+3t\right]$ exactly once.

Consider the differences from the triples of the form $(0$, $a_{i}+p+t$, $b_{i}+p+t)$ with $1 \leq i \leq t$ given by using Construction~\ref{skolem01} with $c=p+t$. By Lemma~\ref{llang1} these triples use edge labels $\left[p+1,p+3t\right]$. If we take the union of the sets of edge labels, we obtain $\left[1,5p+3t\right]$.
%In In a Skolem sequence, the $a_i,b_i,c_i$ and $d_i$ are all unique.

\noindent\textbf{Case 2:} The case $p \equiv 3\ ({\rm mod}\ 4)$, is proved similarly to Case 1. Instead of a Skolem sequence for $S_{1}$, use a hooked Skolem sequence of order $p$ to construct the entries $a_i$ and $b_i$. Instead of a Skolem sequence for $S_{2}$, use a near-Skolem sequence of order $2p+1$ to construct the entries $d_{a_i}$ and $d_{b_i}$. The result follows in a similar fashion.
% of the form $\left(0,a_{i}+p+t,b_{i}+p+t\right)$ with $1 \leq i \leq t$ given by using Construction~\ref{skolem01} ($c=p+1, t=0$). These triples give the vertex labels for $t$ triangles ($C_{3}^{t}$), where $0$ is repeated $t$ times as a common vertex. Thus
%   \textbf{Case 2:} If $p \equiv 3\  (mod\ 4)$, will be the same proof.

For the two remaining cases, we consider $p \equiv 1\ ({\rm mod}\ 4)$ and $p \equiv 2\ ({\rm mod}\ 4)$ where we are using near gracefully labelled $C_{5}^{p}$.

\noindent\textbf{Case 3:} If $p \equiv 1\ ({\rm mod}\ 4)$, use a hooked Skolem sequence of order $p$ for $S_{1}$ to obtain the entries $a_i$ and $b_i$ and use a hooked Skolem sequence of order $2p$ for $S_{2}$ to obtain the entries $d_{a_i}$ and $d_{b_i}$.

\noindent\textbf{Case 4:} If $p \equiv 2\ ({\rm mod}\ 4)$, use a hooked Skolem sequence of order $p$ for $S_{1}$ to obtain the entries $a_i$ and $b_i$ and use a hooked near-Skolem sequence of order $2p+1$ for $S_{2}$ to obtain the entries $d_{a_i}$ and $d_{b_i}$.\end{pf}

We give an example of a labelling of $G=C_{3}^{9}C_{5}^{4}$. Consider the graceful labelling of $C_{5}^{4}$ obtained using Construction~\ref{c5a} with $t=9$ and Construction~\ref{c35} to obtain $\left(0,43,7,8,40\right)$, $\left(0,38,4,2,37\right)$, $\left(0,44,3,6,39\right)$, $\left(0,46,1,5,47\right)$ as the vertex labels of the $C_{5}$ vanes. Apply Construction~\ref{skolem01} with the Langford sequence $\left(13,11,9,7,5,12,10,8,6,5,7,9,11,13,6,8,10,12\right)$ and $c=13$ to construct the triples $\left(0,18,23\right)$, $\left(0,22,28\right)$, $\left(0,17,24\right)$, $\left(0,21,29\right)$, $\left(0,16,25\right)$, $\left(0,20,30\right)$, $\left(0,15,26\right)$, $\left(0,19,31\right)$, $\left(0,14,27\right)$. Together, these give a graceful labelling the graph $G$.

\begin{theorem}\label{thirdthm2}If $G=C_3^{t}C_{5}^{p}$ and $t \geq 2p+1$ then
\begin{enumerate}
	\item $G$ is graceful when $p \equiv 0\ ({\rm mod}\ 4)$ and $t \equiv 0,1 \pmod 4$,
	
	\item $G$ is graceful when $p \equiv 3\ ({\rm mod}\ 4)$ and $t \equiv 0,3 \pmod 4$,
	
	\item $G$ is near graceful when $p \equiv 1\ ({\rm mod}\ 4)$ and $t \equiv 0,3 \pmod 4$,
	
	\item $G$ is near graceful when $p \equiv 2\ ({\rm mod}\ 4)$ and $t \equiv 0,1 \pmod 4$.
	
\end{enumerate}
\end{theorem}
\begin{pf}The (near) graceful $C_{5}^{p}$ exists by Theorem~\ref{thirdthm}. Construction~\ref{c35} gives the edge and vertex labels of $G$. Then $G$ is graceful when $p \equiv 0,3 \pmod 4$ and near graceful when $p \equiv 1,2 \pmod 4$ by Lemma~\ref{c35c}.
\end{pf}
 %From Construction ? we know a Langford sequence of order $t$ and defect $d=p+1$ exists if and only if $t \geq 2d-1$, $t \equiv 0,1 \pmod 4$ and $d$ is odd or $t \equiv 0,3 \pmod 4$ and $d$ is even. Thus, the Langford sequence we need exists then I can use Construction ??.
%\ref{c35n} it has labels as Lemma\ref{c35c},\ref{c35cn}.
Theorem~\ref{thirdthm2} only contains results for half of the possible combinations of $p$ and $t$ and only for large $t$. That is not to say the other cases are not (near) graceful. For example, the case $p=t=1$ does not lend itself to our construction. However the labelling of $C_3^{1}C_{5}^{1}$ with the vertices of the vanes labelled by $\left(0,5,7\right)$ and $\left(0,8,4,3,6\right)$ is graceful.

\section{$C_{3}^{t}C_{6}^{h}$}\label{section6}
%In \cite{dyer}, Dyer, Payne, Shalaby, and Wicks \cite{dyer} used Skolem-type sequences to gracefully label Dutch windmills. Here we will the same method of labellings and derive it to label $C_{3}^{(t)}C_{6}^{(h)}$ windmills.
In this section, we extend the technique of \cite{Bermond} to obtain labellings for $C_{3}^{t}C_{6}^{h}$.
%%%%%%%%%%%%%%%%%%%%%%%%%%%%%%%%%%%%%%%%%%%%%%%%%%%%%%%%%%%%%%%%%%%%%%%%%%%%%%%%%%%%%%%%%%%%%%%%%%%%%%%%%%%%%%%%%%%%%%%%%%%%%%%%%%%%%%%%%%%%%%%%%%%%%%%%%%%%%%%%%%%%%%%%
\begin{lemma}\label{lw36}Suppose there exists a labelling of a windmill with two triangles labelled $\left(0,i,b_{i}+n\right)$ and $\left(0,j,b_{j}+n\right)$, where $\left(a_i,b_i\right)$ and $\left(a_j,b_j\right)$ are the positions of $i$ and $j$ in a Skolem sequence. If these two triangles are removed and replaced by a $C_6$ with vertex labels $\left(0,b_{i}+n,i,i+j,j,b_{j}+n\right)$, for $1 \leq i,j \leq n$ and $i\neq j$, then the edge labels are preserved.\end{lemma}
\begin{pf}The edge labels induced by the triples $\left(0,i,b_{i}+n\right)$ and $\left(0,j,b_{j}+n\right)$ are $\left\{i,j,a_{j}+n,b_{j}+n,a_{i}+n,b_{i}+n\right\}$. The edge labels induced by $\left(0\right.$, $b_{i}+n$, $i$, $i+j$, $j$, $b_{j}+n)$ are the same.\end{pf}
%%%%%%%%%%%%%%%%%%%%%%%%%%%%%%%%%%%%%%%%%%%%%%%%%%%%%%%%%%%%%%Constructions for $C_{3}^{(t)}C_{6}^{(h)}$ %%%%%%%%%%%%%%%%%%%%%%%%%%%%%%%%%%%%%%%%%%%%%%%%%%%%%%%%%%%%%%

\begin{construction}\label{w36}
From a (near) gracefully labelled $C_{3}^{n}$ labelled by $(0$, $i$, $b_{i}+n)$ with $1 \leq i \leq n$ formed by one of the (hooked) Skolem sequences of order $n$ in Table~\ref{table1} - \ref{table4}, in Appendix A as in Construction~\ref{skolem01}(\ref{skolem012}), replace $2h$ triples $\left(h \leq \left\lfloor \frac{(2n+1)}{5}\right\rfloor\right)$ with $h$ 6-tuples by replacing the pair of triples $\left(0,i,b_{i}+n\right)$ and $\left(0,j,b_{j}+n\right)$ with $\left(0,b_{i}+n,i,i+j,j,b_{j}+n\right)$, with pairs indicated as in Table~\ref{ww}.% and this pairing gives us up to $h \leq \left\lfloor \frac{2t+1}{5}\right\rfloor$.
\end{construction}

Table~\ref{ww} gives a family of possible pairs $\left(i,j\right)$ corresponding to triples $\left(0,i,b_{i}+n\right)$ and $\left(0,j,b_{j}+n\right)$ with $1 \leq i,j \leq n$ that may be paired to form hexagons of the form $\left(0,b_{i}+n,i,i+j,j,b_{j}+n\right)$ using the vertex label $i+j$. If the value of $i+j$ does not conflict with any other vertex labels, we obtain a (near) graceful labelling of $C_{3}^{t}C_{6}^{h}$, where $n=t+2h$ and $1 \leq t,h \leq n$.

We notice from Table~\ref{ww}, when $n \in \left\{5k,5k+1\right\}$ with $k \geq 1$, we can obtain up to $2k$ distinct values of $i+j$. It is straightforward to check that this is the best possible result using this method of generating hexagons.

The bound $h \leq \left\lfloor \frac{(2n+1)}{5}\right\rfloor$ comes from the maximum number of possible pairs of triples formed by using the method of Table~\ref{ww}. Combining this with the fact $n=t+2h$, it is straightforward to show these restrictions are equivalent to $h \leq 2t+1$.

\begin{table}[ht]
\begin{center}%\renewcommand{\arraystretch}{1.5}
\scalebox{0.93}{
\begin{tabular}{|c|c|c|}
\hline
$n$    & $\left(i,j\right)$  &  \\ \hline

$5k$   & $\left(k\!+\!z,4k\!+\!z\!+\!1\right)$,\ $\left(2k\!+\!z'\!+\!1,3k\!+\!z'\!+\!1\right)$ & $0 \leq z,z' \leq k\!-\!1$ \\ \hline

$5k\!+\!1$ & $\left(k\!+\!z\!+\!1,4k\!+\!z\!+\!1\right)$,\ $\left(2k\!+\!z'\!+\!2,3k\!+\!z'\!+\!1\right)$ & $0 \leq z \leq k$,\ $0 \leq z' \leq k\!-\!2$ \\ \hline

$5k\!+\!2$ & $\left(k\!+\!z\!+\!1,4k\!+\!z\!+\!2\right)$,\ $\left(2k\!+\!z'\!+\!2,3k\!+\!z'\!+\!2\right)$ & $0 \leq z \leq k$,\ $0 \leq z' \leq k\!-\!1$ \\ \hline

$5k\!+\!3$ & $\left(k\!+\!z\!+\!1,4k\!+\!z\!+\!3\right)$,\ $\left(2k\!+\!z'\!+\!2,3k\!+\!z'\!+\!3\right)$ & $0 \leq z \leq k$,\ $0 \leq z' \leq k\!-\!1$ \\ \hline

$5k\!+\!4$ & $\left(k\!+\!z\!+\!1,4k\!+\!z\!+\!4\right)$,\ $\left(2k\!+\!z'\!+\!3,3k\!+\!z'\!+\!3\right)$ & $0 \leq z \leq k$,\ $0 \leq z' \leq k\!-\!1$ \\ \hline

\end{tabular}}
\caption{Possible pairs of triples corresponding to $i$ and $j$ that may be paired to form hexagons with $i+j$.} \label{ww}
\end{center}
\end{table}

Recall from the example following Lemma~\ref{llang1} the hooked Skolem sequence yielded the triples $\left(0,1,6\right),\left(0,2,10\right)$, and $\left(0,3,7\right)$. Pair the triples $\left(0,1,6\right)$ and $\left(0,3,7\right)$ to form the $6$-tuple $\left(0,6,1,4,3,7\right)$ using the label $4$ which was not already used, so we obtain a near graceful labelling of $C_{3}^{1}C_{6}^{1}$. We cannot pair $\left(0,1,6\right)$ with $\left(0,2,10\right)$ because $1+2=3$ which duplicates another vertex label.

We now present how to use Construction~\ref{w36} with Table~\ref{ww} to get a family of possible pairs of triples corresponding to $i$ and $j$ that may be paired to form hexagons with $i+j$. Start with the Skolem sequence $(8$, $6$, $4$, $2$, $7$, $2$, $4$, $6$, $8$, $3$, $5$, $7$, $3$, $1$, $1$, $5)$ constructed from Table~\ref{table1} in Appendix A. Take the triples of the form $\left(0,i,b_{i}+8\right)$. From Table~\ref{ww}, we can get up to three pairs $\left(2,7\right),\left(3,8\right)$, and $\left(4,6\right)$ which give three different gracefully labelled graphs. We have three possible replacements: replace $\left(0,2,14\right)$ and $\left(0,7,20\right)$ by $\left(0,14,2,9,7,20\right)$; replace $\left(0,3,21\right)$ and $\left(0,8,17\right)$ by $\left(0,21,3,11,8,17\right)$; and replace $\left(0,4,15\right)$ and $\left(0,6,16\right)$ by $\left(0,15,4,10,6,16\right)$.
We can gracefully label $C_{3}^{6}C_{6}^{1}$ by using any one of these replacements, $C_{3}^{4}C_{6}^{2}$ by using any two, and $C_{3}^{2}C_{6}^{3}$ by simultaneously using all three.

Note that the first right endpoint in the Skolem sequences of order $n$ described in Tables~\ref{table1}--\ref{table4}, in Appendix A, is always at $\left\lceil \frac{(n+3)}{2} \right\rceil$. This fact will be useful later in the proof of Theorem~\ref{tw63}.

In Theorem~\ref{tw63}, we use Construction~\ref{w36} with the appropriate Skolem sequence of order $n$ which will give (near) graceful labellings for $C_{3}^{t}C_{6}^{h}$ and $h \leq 2t+1$ as given in Table~\ref{c36r}.
\begin{table}[ht!]
\begin{center}%\renewcommand{\arraystretch}{1.5}
\scalebox{1.00}{
\begin{tabular}{|c|c|c|}
\hline
$n=t+2h$                 & Graceful                               & Near graceful \\ \hline
$5k$   &$k \equiv 0,1\ ({\rm mod}\ 4)$          &$k \equiv 2,3\ ({\rm mod}\ 4)$               \\ \hline
$5k+1$ &$k \equiv 0,3\ ({\rm mod}\ 4)$          &$k \equiv 1,2\ ({\rm mod}\ 4)$               \\ \hline
$5k+2$ &$k \equiv 2,3\ ({\rm mod}\ 4)$          &$k \equiv 0,1\ ({\rm mod}\ 4)$               \\ \hline
$5k+3$ &$k \equiv 1,2\ ({\rm mod}\ 4)$          &$k \equiv 0,3\ ({\rm mod}\ 4)$               \\ \hline
$5k+4$ &$k \equiv 0,1\ ({\rm mod}\ 4)$          &$k \equiv 2,3\ ({\rm mod}\ 4)$               \\ \hline
\end{tabular}}
\caption{Summary of results of (near) graceful labelling $C_{3}^{t}C_{6}^{h}$ and $h \leq 2t+1$.}\label{c36r}
\end{center}
\end{table}
\begin{theorem}\label{tw63}If $G=C_{3}^{t}C_{6}^{h}$ with $h \leq 2t+1$, then $G$ is graceful or near graceful.\end{theorem}

\begin{pf}Let $n=2h+t$. We begin by considering $n=5k$. Proofs for other $n$ follow in the same fashion.

Start with a (near) gracefully labelled $C_{3}^{n}$ with $n=5k$ with $k \equiv 0,1\!\!\pmod{4}$ labelled by $\left(0,i,b_{i}+n\right)$ with $1 \leq i \leq n$ formed by a Skolem sequence as indicated in Construction~\ref{w36}, where $2h+t=n$ and $h \leq 2t+1$. By Construction~\ref{w36}, we replace $2h$ triples with $h$ 6-tuples by $(0$, $b_{i}+n$, $i$, $i+j$, $j$, $b_{j}+n)$. By Lemma~\ref{lw36} the edge labels will not be changed. Therefore, in this new labelling we use the edge labels $\left[1,3n\right]$. The new labelling uses the same vertex labels, as well as some new labels of the form $i+j$.

We constructed $\left(0,i,b_{i}+n\right)$ with $1 \leq i \leq n$ by using Table~\ref{table1} or Table~\ref{table2}, in Appendix A. These triples use the elements $i$ in the interval $\left[1,n\right]$ and use the elements $b_{i}+n$ in the interval $\left[\left\lceil \frac{(3n+3)}{2} \right\rceil,3n\right]$. Thus any subset of these triples do not use any vertex labels in the interval $\left[n+1,\left\lceil \frac{(3n+3)}{2} \right\rceil -1 \right]$.
% Since is the first right end point at position $b_{2}=\left\lceil \frac{n+3}{2} n \right\rceil$ and $n=5k$ with $i \leq 5k$ we did not use any labels in the interval $\left[5k+1,7k+1\right]$. (expanded) middle entry first, last entry, range.
% From Table~\ref{ww}, we know all these vertices are distinct and in the interval $\left[t+1,\left\lceil \frac{3t+3}{2} \right\rceil -1 \right]$.

By construction, the sums $i+j$ are all distinct and in the interval\linebreak $\left[n+1,\left\lceil \frac{(3n+3)}{2} \right\rceil -1 \right]$. The new vertex labels do not duplicate any previously used labels; this new vertex labelling is injective and uses labels from the set $\left[1,3n\right]$.

Finally, we can conclude that since the vertex labels are a subset of $\left[0,3n\right]$ and the edge labels are exactly $\left[1,3n\right]$, that $G=C_{3}^{t}C_{6}^{h}$ with $2h+t=n$ and $h \leq 2t+1$ can be gracefully labelled when $k \equiv 0,1\ ({\rm mod}\ 4)$.
%$\left[5k+1,7k+1\right]$
%We can get $2k$ distinct value of $i+j$, by the pairing method in the table we got this bound $\left\lfloor \frac{2n+1}{5}\right\rfloor$ and can't work if $h$, greater than this bound.
%%%%%%%%%%%%%%%%%%%%%%%%%%%%%%%%%%%%%%%%%%%%%%%%%

The case $n=5k$ with $k \equiv 2,3\ ({\rm mod}\ 4)$, works in the same way except the construction of $\left(0,i,b_{i}+n\right)$ with $1 \leq i \leq n$ is by using Table~\ref{table3} or Table~\ref{table4}, in Appendix A. This construction then uses the edge labels $\left[1,3n-1\right] \cup \left\{3n+1\right\}$ and vertex labels from $\left[0,3n-1\right]  \cup \left\{3n+1\right\}$. We conclude that $G=C_{3}^{t}C_{6}^{h}$ with $2h+t=5k$ and $h \leq 2t+1$ can be near gracefully labelled when $k \equiv 2,3\ ({\rm mod}\ 4)$.\end{pf}
Theorem~\ref{tw63} only contains results for $h \leq 2t+1$. We cannot (near) gracefully label $G=C_{3}^{t}C_{6}^{h}$ with $h > 2t+1$ by our construction, so these problems remain open.

\section{Discussion}\label{section7}

In this paper, we have completely characterized the situation where $C_{3}^{t}C_{4}^{s}$ is graceful or near graceful, and given partial solutions for $C_{3}^{t}C_{5}^{p}$ and $C_{3}^{t}C_{6}^{h}$.

Many of the techniques of this paper can be combined, but are difficult to reduce to theorems. From the example in Section~\ref{section5} we gracefully labelled $C_{3}^{9}C_{5}^{4}$. If we replace the triples $\left(0,16,25\right),\left(0,20,30\right), \left(0,19,31\right)$ by $\left(0,9,25\right),\left(0,10,30\right),\left(0,12,31\right)$, we notice that the edge labels are the same and the new vertex labels do not appear elsewhere in the labelling of $C_{3}^{9}C_{5}^{4}$.  Thus we have obtained another graceful labelling of $C_{3}^{9}C_{5}^{4}$. Now by the same technique we used in Section~\ref{section6}, we can obtain a graceful labelling for $C_{3}^{7} C_{5}^{4} C_{6}^{1}$ by replacing $\left(0,9,25\right),\left(0,10,30\right)$ with $\left(0,25,9,19,10,30\right)$.

In \cite{yang7,yang9,yang11,yang13}, it was proved that graceful labellings exist for $C_{7}^{e}$, $C_{9}^{e}$, $C_{11}^{e}$, and $C_{13}^{e}$ respectively. By Theorem 4.3(2) in \cite{ahmad}, we can prove that graceful labellings exist for the $C_{3}^{t}C_{7}^{e}$, $C_{3}^{t}C_{9}^{e}$, $C_{3}^{t}C_{11}^{e},$ and $C_{3}^{t}C_{13}^{e}$, respectively with $t$ sufficiently larger than $e$. It remains to determine the gracefulness of these windmills when $t$ is small.
%We should include an expand ideas (take an exists results and see if we can add our labelling, for example, 3-11 windmill, like that theorem with Nabil paper)

Since Skolem-type sequences have proven useful in (near) gracefully labelling $C_3$ windmills, as well as enabling us to label $C_3C_4$ windmills, it is natural to ask can we use Skolem-type sequences to gracefully label all $C_{4}$ windmills? While the result is known, this would be an extension of the Skolem techniques used in this paper.

In the process of writing this paper, many Skolem-type sequences were considered. We pose two open problems whose solutions would be of interest not only to these constructions, but in their own right. First, can we find a family of $m$-fold Skolem-type sequences with first right endpoint as large as possible? Second, what are necessary and sufficient conditions for the existence of $m$-fold Langford sequences with $m \geq 2$?

\section*{Appendix A: Skolem Constructions}

\begin{table}[H]
\begin{center}%\setlength\extrarowheight{7pt}
\scalebox{1.00}{
\begin{tabular}{|c|c|c|c|c|}
\hline
  $i$ & $a_i$ & $b_i$ &  \\
\hline
\hline
 $2r+2$ &\ $2m-r$   & $2m+2+r$ & $0 \leq r \leq 2m-1$  \\
\hline
 $1$ &\ $7m$   & $7m+1$ &  $-$  \\
\hline
 $4m-1$ &\ $2m+1$   & $6m$ &  $-$  \\
\hline
 $2m+2r+1$ &\ $5m+1-r$   & $7m+r+2$ &  $0 \leq r \leq m-2$  \\
\hline
 $2m-1$ &\ $4m+2$   & $6m+1$ &  $-$  \\
\hline
 $2m-3-2r$ &\ $5m+2+r$   & $7m-1-r$ &  $0 \leq r \leq m-3$  \\
\hline
\end{tabular}}
\caption{Skolem sequence construction of order $n=4m$ and $m\geq 1$ from \cite{skolem}.} \label{table1}
\end{center}
\end{table} \vspace{-3ex}
\begin{table}[H]
\begin{center}%\setlength\extrarowheight{7pt}
\scalebox{1.00}{
\begin{tabular}{|c|c|c|c|c|}
\hline
 $i$ & $a_i$ & $b_i$ &  \\
\hline
\hline
$2r$ &\ $2m+1-r$   & $2m+1+r$ & $1 \leq r \leq 2m$  \\
\hline
 $4m+1$ &\ $2m+1$   & $6m+2$ &  $-$  \\
\hline
 $2m-1+2r$ &\ $5m+2-r$   & $7m+1+r$ &  $1 \leq r \leq m$  \\
\hline
 $2m-1$ &\ $6m+3$   & $8m+2$ &  $-$  \\
\hline
 $1$ &\ $5m+2$   & $5m+3$ &  $-$  \\
\hline
 $2r+1$ &\ $6m+2-r$   & $6m+3+r$ &  $1 \leq r \leq m-2$  \\
\hline
\end{tabular}}
\caption{Skolem sequence construction of order $n=4m+1$and $m \geq 2$ from \cite{marsh}.} \label{table2}
\end{center}
\end{table} \vspace{-3ex}
\begin{table}[H]
\begin{center}%\setlength\extrarowheight{7pt}
\scalebox{1.00}{
\begin{tabular}{|c|c|c|c|c|}
\hline
 $i$ & $a_i$ & $b_i$ &  \\
\hline
\hline
 $2r$ &\ $2m+2-r$   & $2m+2+r$ & $1 \leq r \leq 2m+1$  \\
\hline
 $1$ &\ $7m+4$   & $7m+5$ &  $-$  \\
\hline
 $1+2r$ &\ $6m+2-r$   & $6m+3+r$ &  $1 \leq r \leq m$  \\
\hline
 $2m+3$ &\ $6m+2$   & $8m+5$ &  $-$  \\
\hline
 $2m+3+2r$ &\ $5m+2-r$   & $7m+5+r$ &  $1 \leq r \leq m-2$  \\
\hline
 $4m+1$ &\ $2m+2$   & $6m+3$ &  $-$  \\
\hline
\end{tabular}}
\caption{Hooked Skolem sequence construction of order $n=4m+2$and $m \geq 2$ from \cite{hilton}.} \label{table3}
\end{center}
\end{table} \vspace{-3ex}
\begin{table}[H]
\begin{center}%\setlength\extrarowheight{7pt}
\scalebox{1.00}{
\begin{tabular}{|c|c|c|c|c|}
\hline
 $i$ & $a_i$ & $b_i$ &  \\
\hline
\hline
 $2r$ &\ $2m+2-r$   & $2m+2+r$ & $1 \leq r \leq 2m+1$  \\
\hline
 $1$ &\ $5m+4$   & $5m+5$ &  $-$  \\
\hline
 $1+2r$ &\ $6m+5-r$   & $6m+6+r$ &  $1 \leq r \leq m-1$  \\
\hline
 $2m+1$ &\ $6m+6$   & $8m+7$ &  $-$  \\
\hline
 $2m+1+2r$ &\ $5m+4-r$   & $7m+5+r$ &  $1 \leq r \leq m$  \\
\hline
 $4m+3$ &\ $2m+2$   & $6m+5$ &  $-$  \\
\hline
\end{tabular}}
\caption{Hooked Skolem sequence construction of order $n=4m+3$ and $m \geq 1$ from \cite{hilton}.} \label{table4}
\end{center}
\end{table}
%\end{appendices}
%%%%%%%%%%%%%%%%%%%%%%%%%%%%%%%%%%%%%%%%%%%%%%%%%%%%%%%%%%%%%%%%%%%xxxxxxxxxxxxxx%%%%%%%%%%%%%%%%%%%%%%%%%%%%%%%%%%%%%%%%%%%%%%%%%%%%%%%%%%%%%%%%%%%%%%%%%%%%%%%%%%%%%%%
%\section*{Appendix B}
\section*{Appendix B}
The following table includes graceful labellings of $C_{3}^{1}C_{4}^{s}$, where $1 \leq s \leq 20$.

The vanes $(0,3,2,6)$ and $(0,5,7)$ are included with each of the sets in the following table. Note that these two vanes give a graceful labelling when $s=1$.

\begin{longtable}{|l|>{\raggedright\arraybackslash}p{0.9\linewidth}|}
\hline
$s$ & Additional vanes \\
\hline
$2$ & $\left(0,9,1,11\right)$ \\
\hline
$3$ & $\left(0,11,1,15\right)$, $\left(0,12,4,13\right)$ \\
\hline
$4$ & $\left(0,13,1,19\right)$, $\left(0,14,4,15\right)$, $\left(0,16,8,17\right)$ \\
\hline
$5$ & $\left(0,15,1,19\right)$, $\left(0,16,4,17\right)$, $\left(0,20,11,22\right)$, $\left(0,21,13,23\right)$ \\
\hline
$6$ & $\left(0,17,1,21\right)$, $\left(0,18,4,19\right)$, $\left(0,22,12,25\right)$, $\left(0,23,14,26\right)$, $\left(0,24,16,27\right)$ \\
\hline
$7$ & $\left(0,19,1,21\right)$, $\left(0,22,10,27\right)$, $\left(0,23,12,28\right)$, $\left(0,24,14,29\right)$, $\left(0,25,16,30\right)$, $\left(0,26,18,31\right)$ \\
\hline
$8$ & $\left(0,23,1,39\right)$, $\left(0,24,10,31\right)$, $\left(0,25,12,32\right)$, $\left(0,26,14,33\right)$, $\left(0,27,16,34\right)$, $\left(0,28,18,35\right)$, $\left(0,29,20,36\right)$, $\left(0,30,22,37\right)$ \\
\hline
$9$ & $\left(0,23,1,39\right)$, $\left(0,24,10,31\right)$, $\left(0,25,12,32\right)$, $\left(0,26,14,33\right)$, $\left(0,27,16,34\right)$, $\left(0,28,18,35\right)$, $\left(0,29,20,36\right)$, $\left(0,30,22,37\right)$ \\
\hline
$10$ & $\left(0,23,1,41\right)$, $\left(0,42,4,43\right)$, $\left(0,24,10,31\right)$, $\left(0,25,12,32\right)$, $\left(0,26,14,33\right)$, $\left(0,27,16,34\right)$, $\left(0,28,18,35\right)$, $\left(0,29,20,36\right)$, $\left(0,30,22,37\right)$, $\left(0,3,2,6\right)$. $\left(0,5,7\right)$\\
\hline
$11$ & $\left(0,28,13,36\right)$, $\left(0,29,15,37\right)$, $\left(0,30,17,38\right)$, $\left(0,44,18,45\right)$, $\left(0,31,19,39\right)$, $\left(0,32,21,40\right)$, $\left(0,46,22,47\right)$, $\left(0,33,23,41\right)$, $\left(0,34,25,42\right)$, $\left(0,35,27,43\right)$\\
\hline
$12$ & $\left(0,34,26,39\right)$, $\left(0,33,24,38\right)$, $\left(0,32,22,37\right)$, $\left(0,31,20,36\right)$, $\left(0,30,18,35\right)$, $\left(0,45,27,51\right)$, $\left(0,44,25,50\right)$, $\left(0,43,23,49\right)$, $\left(0,42,21,48\right)$, $\left(0,41,19,47\right)$, $\left(0,40,17,46\right)$ \\
\hline
$13$ & $\left(0,38,26,43\right)$, $\left(0,37,24,42\right)$, $\left(0,36,22,41\right)$, $\left(0,35,20,40\right)$, $\left(0,34,18,39\right)$, $\left(0,49,27,55\right)$, $\left(0,48,25,54\right)$, $\left(0,47,23,53\right)$, $\left(0,46,21,52\right)$, $\left(0,45,19,51\right)$, $\left(0,44,17,50\right)$\\
\hline
$14$ & $\left(0,42,26,47\right)$, $\left(0,41,24,46\right)$, $\left(0,40,22,45\right)$, $\left(0,39,20,44\right)$, $\left(0,38,18,43\right)$, $\left(0,53,27,59\right)$, $\left(0,52,25,58\right)$, $\left(0,51,23,57\right)$, $\left(0,50,21,56\right)$, $\left(0,49,19,55\right)$, $\left(0,48,17,54\right)$\\
\hline
$15$ & $\left(0,55,1,59\right)$, $\left(0,56,4,57\right)$, $\left(0,60,10,61\right)$, $\left(0,28,13,36\right)$, $\left(0,62,14,63\right)$, $\left(0,29,15,37\right)$, $\left(0,30,17,38\right)$, $\left(0,44,18,45\right)$, $\left(0,31,19,39\right)$, $\left(0,32,21,40\right)$, $\left(0,46,22,47\right)$, $\left(0,33,23,41\right)$, $\left(0,34,25,42\right)$, $\left(0,35,27,43\right)$\\
\hline
$16$ & $\left(0,44,36,51\right)$, $\left(0,43,34,50\right)$, $\left(0,42,32,49\right)$, $\left(0,41,30,48\right)$, $\left(0,40,28,47\right)$, $\left(0,39,26,46\right)$, $\left(0,38,24,45\right)$, $\left(0,59,37,67\right)$, $\left(0,58,35,66\right)$, $\left(0,57,33,65\right)$, $\left(0,56,31,64\right)$, $\left(0,55,29,63\right)$, $\left(0,54,27,62\right)$, $\left(0,53,25,61\right)$, $\left(0,52,23,60\right)$\\
\hline
$17$ &  $\left(0,48,36,55\right)$, $\left(0,47,34,54\right)$, $\left(0,46,32,53\right)$, $\left(0,45,30,52\right)$, $\left(0,44,28,51\right)$, $\left(0,43,26,50\right)$, $\left(0,42,24,49\right)$, $\left(0,63,37,71\right)$, $\left(0,62,35,70\right)$, $\left(0,61,33,69\right)$, $\left(0,60,31,68\right)$, $\left(0,59,29,67\right)$, $\left(0,58,27,66\right)$, $\left(0,57,25,65\right)$, $\left(0,56,23,64\right)$\\
\hline
$18$ & $\left(0,52,36,59\right)$, $\left(0,51,34,58\right)$, $\left(0,50,32,57\right)$, $\left(0,49,30,56\right)$, $\left(0,48,28,55\right)$, $\left(0,47,26,54\right)$, $\left(0,46,24,53\right)$, $\left(0,67,37,75\right)$, $\left(0,66,35,74\right)$, $\left(0,65,33,73\right)$, $\left(0,64,31,72\right)$, $\left(0,63,29,71\right)$, $\left(0,62,27,70\right)$, $\left(0,61,25,69\right)$, $\left(0,60,23,68\right)$\\
\hline
$19$ & $\left(0,56,36,63\right)$, $\left(0,55,34,62\right)$, $\left(0,54,32,61\right)$, $\left(0,53,30,60\right)$, $\left(0,52,28,59\right)$, $\left(0,51,26,58\right)$, $\left(0,50,24,57\right)$, $\left(0,71,37,79\right)$, $\left(0,70,35,78\right)$, $\left(0,69,33,77\right)$, $\left(0,68,31,76\right)$, $\left(0,67,29,75\right)$, $\left(0,66,27,74\right)$, $\left(0,65,25,73\right)$, $\left(0,64,23,72\right)$\\
\hline
$20$ & $\left(0,61,1,80\right)$, $\left(0,81,4,82\right)$, $\left(0,73,9,83\right)$, $\left(0,65,10,69\right)$, $\left(0,66,12,70\right)$, $\left(0,75,13,76\right)$, $\left(0,67,14,71\right)$, $\left(0,68,16,72\right)$, $\left(0,40,17,46\right)$, $\left(0,30,18,35\right)$, $\left(0,41,19,47\right)$, $\left(0,31,20,36\right)$, $\left(0,42,21,48\right)$, $\left(0,32,22,37\right)$, $\left(0,43,23,49\right)$, $\left(0,33,24,38\right)$, $\left(0,44,25,50\right)$, $\left(0,34,26,39\right)$, $\left(0,45,27,51\right)$\\
\hline	
\end{longtable}

The following table includes near graceful labellings of $C_{3}^{2}C_{4}^{s}$, where $1 \leq s \leq 20$.

\begin{longtable}{|l|>{\raggedright\arraybackslash}p{0.9\linewidth}|}
\hline
$s$ & Vanes \\
\hline
$1$ & $\left(0,5,2,6\right)$, $\left(0,7,8\right)$, $\left(0,9,11\right)$ \\
\hline
$2$ & $\left(0,10,1,12\right)$, $\left(0,5,2,6\right)$, $\left(0,7,8\right)$, $\left(0,15,13\right)$\\
\hline
$3$ & $\left(0,12,1,16\right)$, $\left(0,13,4,14\right)$, $\left(0,5,2,6\right)$,$\left(0,7,8\right)$, $\left(0,19,17\right)$\\
\hline
$4$ & $\left(0,14,1,20\right)$, $\left(0,15,4,16\right)$, $\left(0,17,8,18\right)$, $\left(0,5,2,6\right)$, $\left(0,7,8\right)$, $\left(0,23,21\right)$\\
\hline
$5$ & $\left(0,16,1,20\right)$, $\left(0,17,4,18\right)$, $\left(0,21,11,23\right)$, $\left(0,22,13,24\right)$, $\left(0,5,2,6\right)$, $\left(0,7,8\right)$, $\left(0,27,25\right)$\\
\hline
$6$ & $\left(0,18,1,22\right)$, $\left(0,19,4,20\right)$, $\left(0,23,12,26\right)$, $\left(0,24,14,27\right)$,$\left(0,25,16,28\right)$, $\left(0,5,2,6\right)$, $\left(0,7,8\right)$, $\left(0,31,29\right)$\\
\hline
$7$ & $\left(0,20,1,22\right)$, $\left(0,23,10,28\right)$, $\left(0,24,12,29\right)$, $\left(0,25,14,30\right)$, $\left(0,26,16,31\right)$, $\left(0,27,18,32\right)$, $\left(0,5,2,6\right)$, $\left(0,7,8\right)$, $\left(0,33,35\right)$\\
\hline
$8$ & $\left(0,25,16,28\right)$, $\left(0,24,14,27\right)$, $\left(0,23,12,26\right)$, $\left(0,32,17,36\right)$, $\left(0,31,15,35\right)$, $\left(0,30,13,34\right)$, $\left(0,29,11,33\right)$, $\left(0,5,2,6\right)$, $\left(0,7,8\right)$, $\left(0,39,37\right)$\\
\hline
$9$ & $\left(0,24,1,40\right)$, $\left(0,25,10,32\right)$, $\left(0,26,12,33\right)$, $\left(0,27,14,34\right)$, $\left(0,28,16,35\right)$, $\left(0,29,18,36\right)$, $\left(0,30,20,37\right)$, $\left(0,31,22,38\right)$, $\left(0,5,2,6\right)$, $\left(0,7,8\right)$, $\left(0,43,41\right)$\\
\hline
$10$ & $\left(0,24,1,42\right)$, $\left(0,43,4,44\right)$, $\left(0,25,10,32\right)$, $\left(0,26,12,33\right)$, $\left(0,27,14,34\right)$, $\left(0,28,16,35\right)$, $\left(0,29,18,36\right)$, $\left(0,30,20,37\right)$, $\left(0,31,22,38\right)$, $\left(0,5,2,6\right)$, $\left(0,7,8\right)$,$\left(0,47,45\right)$\\
\hline
$11$ & $\left(0,29,13,37\right)$, $\left(0,30,15,38\right)$, $\left(0,31,17,39\right)$, $\left(0,45,18,46\right)$, $\left(0,32,19,40\right)$, $\left(0,33,21,41\right)$, $\left(0,47,22,48\right)$, $\left(0,34,23,42\right)$, $\left(0,35,25,43\right)$, $\left(0,36,27,44\right)$, $\left(0,5,2,6\right)$, $\left(0,7,8\right)$, $\left(0,51,49\right)$\\
\hline
$12$ & $\left(0,35,26,40\right)$, $\left(0,34,24,39\right)$, $\left(0,33,22,38\right)$, $\left(0,32,20,37\right)$, $\left(0,31,18,36\right)$, $\left(0,46,27,52\right)$, $\left(0,45,25,51\right)$, $\left(0,44,23,50\right)$, $\left(0,43,21,49\right)$, $\left(0,42,19,48\right)$, $\left(0,41,17,47\right)$, $\left(0,5,2,6\right)$, $\left(0,7,8\right)$, $\left(0,55,53\right)$\\
\hline
$13$ &  $\left(0,39,26,44\right)$, $\left(0,38,24,43\right)$, $\left(0,37,22,42\right)$, $\left(0,36,20,41\right)$, $\left(0,35,18,40\right)$, $\left(0,50,27,56\right)$, $\left(0,49,25,55\right)$, $\left(0,48,23,54\right)$, $\left(0,47,21,53\right)$, $\left(0,46,19,52\right)$, $\left(0,45,17,51\right)$, $\left(0,10,1,12\right)$, $\left(0,5,2,6\right)$, $\left(0,7,8\right)$, $\left(0,59,57\right)$\\
\hline
$14$ & $\left(0,43,26,48\right)$, $\left(0,42,24,47\right)$, $\left(0,41,22,46\right)$, $\left(0,40,20,45\right)$, $\left(0,39,18,44\right)$, $\left(0,54,27,60\right)$, $\left(0,53,25,59\right)$, $\left(0,52,23,58\right)$, $\left(0,51,21,57\right)$, $\left(0,50,19,56\right)$, $\left(0,49,17,55\right)$, $\left(0,12,1,16\right)$, $\left(0,13,4,14\right)$, $\left(0,5,2,6\right)$, $\left(0,7,8\right)$, $\left(0,63,61\right)$\\
\hline
$15$ & $\left(0,56,1,60\right)$, $\left(0,57,4,58\right)$, $\left(0,61,10,62\right)$, $\left(0,29,13,37\right)$, $\left(0,63,14,64\right)$, $\left(0,30,15,38\right)$, $\left(0,31,17,39\right)$, $\left(0,45,18,46\right)$, $\left(0,32,19,40\right)$, $\left(0,33,21,41\right)$, $\left(0,47,22,48\right)$, $\left(0,34,23,42\right)$, $\left(0,35,25,43\right)$, $\left(0,36,27,44\right)$, $\left(0,5,2,6\right)$, $\left(0,7,8\right)$, $\left(0,67,65\right)$\\
\hline
$16$ & $\left(0,45,36,52\right)$, $\left(0,44,34,51\right)$, $\left(0,43,32,50\right)$, $\left(0,42,30,49\right)$, $\left(0,41,28,48\right)$, $\left(0,40,26,47\right)$, $\left(0,39,24,46\right)$, $\left(0,60,37,68\right)$, $\left(0,59,35,67\right)$, $\left(0,58,33,66\right)$, $\left(0,57,31,65\right)$, $\left(0,56,29,64\right)$, $\left(0,55,27,63\right)$, $\left(0,54,25,62\right)$, $\left(0,53,23,61\right)$, $\left(0,5,2,6\right)$, $\left(0,7,8\right)$, $\left(0,71,69\right)$\\
\hline
$17$ & $\left(0,49,36,56\right)$, $\left(0,48,34,55\right)$, $\left(0,47,32,54\right)$, $\left(0,46,30,53\right)$, $\left(0,45,28,52\right)$, $\left(0,44,26,51\right)$, $\left(0,43,24,50\right)$, $\left(0,64,37,72\right)$, $\left(0,63,35,71\right)$, $\left(0,62,33,70\right)$, $\left(0,61,31,69\right)$, $\left(0,60,29,68\right)$, $\left(0,59,27,67\right)$, $\left(0,58,25,66\right)$, $\left(0,57,23,65\right)$, $\left(0,10,1,12\right)$, $\left(0,5,2,6\right)$, $\left(0,7,8\right)$, $\left(0,75,73\right)$\\
\hline
$18$ & $\left(0,53,36,60\right)$, $\left(0,52,34,59\right)$, $\left(0,51,32,58\right)$, $\left(0,50,30,57\right)$, $\left(0,49,28,56\right)$, $\left(0,48,26,55\right)$, $\left(0,47,24,54\right)$, $\left(0,68,37,76\right)$, $\left(0,67,35,75\right)$, $\left(0,66,33,74\right)$, $\left(0,65,31,73\right)$, $\left(0,64,29,72\right)$, $\left(0,63,27,71\right)$, $\left(0,62,25,70\right)$, $\left(0,61,23,69\right)$, $\left(0,12,1,16\right)$, $\left(0,13,4,14\right)$, $\left(0,5,2,6\right)$, $\left(0,7,8\right)$, $\left(0,79,77\right)$\\
\hline
$19$ & $\left(0,57,36,64\right)$, $\left(0,56,34,63\right)$, $\left(0,55,32,62\right)$, $\left(0,54,30,61\right)$, $\left(0,53,28,60\right)$, $\left(0,52,26,59\right)$, $\left(0,51,24,58\right)$, $\left(0,72,37,80\right)$, $\left(0,71,35,79\right)$, $\left(0,70,33,78\right)$, $\left(0,69,31,77\right)$, $\left(0,68,29,76\right)$, $\left(0,67,27,75\right)$, $\left(0,66,25,74\right)$, $\left(0,65,23,73\right)$, $\left(0,14,1,20\right)$, $\left(0,15,4,16\right)$, $\left(0,17,8,18\right)$, $\left(0,5,2,6\right)$, $\left(0,7,8\right)$, $\left(0,83,81\right)$\\
\hline
$20$ & $\left(0,62,1,81\right)$, $\left(0,82,4,83\right)$, $\left(0,74,9,84\right)$, $\left(0,66,10,70\right)$, $\left(0,67,12,71\right)$, $\left(0,76,13,77\right)$, $\left(0,68,14,72\right)$, $\left(0,69,16,73\right)$, $\left(0,41,17,47\right)$, $\left(0,31,18,36\right)$, $\left(0,42,19,48\right)$, $\left(0,32,20,37\right)$, $\left(0,43,21,49\right)$, $\left(0,33,22,38\right)$, $\left(0,44,23,50\right)$, $\left(0,34,24,39\right)$, $\left(0,45,25,51\right)$, $\left(0,35,26,40\right)$, $\left(0,46,27,52\right)$, $\left(0,5,2,6\right)$, $\left(0,7,8\right)$, $\left(0,87,85\right)$\\
\hline
\end{longtable}

In the following table we include near graceful labellings of $C_{3}^{3}C_{4}^{s}$, where $1 \leq s \leq 19$.

\begin{longtable}{|l|>{\raggedright\arraybackslash}p{0.9\linewidth}|}
\hline
$s$ & Vanes \\
\hline
$1$ & $\left(0,8,2,9\right)$, $\left(0,3,5\right)$, $\left(0,11,12\right)$, $\left(0,10,14\right)$\\
\hline
$2$ & $\left(0,11,1,13\right)$, $\left(0,8,2,9\right)$, $\left(0,3,5\right)$, $\left(0,14,18\right)$, $\left(0,15,16\right)$\\
\hline
$3$ & $\left(0,13,1,17\right)$, $\left(0,14,4,15\right)$, $\left(0,8,2,9\right)$, $\left(0,3,5\right)$, $\left(0,22,18\right)$, $\left(0,20,19\right)$\\
\hline
$4$ & $\left(0,15,1,19\right)$, $\left(0,16,4,17\right)$, $\left(0,20,10,21\right)$, $\left(0,8,2,9\right)$, $\left(0,3,5\right)$, $\left(0,26,22\right)$, $\left(0,24,23\right)$\\
\hline
$5$ & $\left(0,17,1,21\right)$, $\left(0,18,4,19\right)$, $\left(0,22,11,24\right)$, $\left(0,23,13,25\right)$, $\left(0,8,2,9\right)$, $\left(0,3,5\right)$, $\left(0,30,26\right)$, $\left(0,28,27\right)$\\
\hline
$6$ & $\left(0,19,1,23\right)$, $\left(0,20,4,21\right)$, $\left(0,24,12,27\right)$, $\left(0,25,14,28\right)$, $\left(0,26,16,29\right)$, $\left(0,8,2,9\right)$, $\left(0,3,5\right)$, $\left(0,34,30\right)$, $\left(0,32,31\right)$\\
\hline
$7$ & $\left(0,21,1,23\right)$, $\left(0,24,10,29\right)$, $\left(0,25,12,30\right)$,$\left(0,26,14,31\right)$, $\left(0,27,16,32\right)$, $\left(0,28,18,33\right)$, $\left(0,8,2,9\right)$, $\left(0,3,5\right)$, $\left(0,38,34\right)$, $\left(0,36,35\right)$\\
\hline
$8$ & $\left(0,26,16,29\right)$, $\left(0,25,14,28\right)$, $\left(0,24,12,27\right)$, $\left(0,33,17,37\right)$, $\left(0,32,15,36\right)$, $\left(0,31,13,35\right)$, $\left(0,30,11,34\right)$, $\left(0,8,2,9\right)$, $\left(0,3,5\right)$, $\left(0,42,38\right)$, $\left(0,40,39\right)$\\
\hline
$9$ & $\left(0,25,1,41\right)$, $\left(0,26,10,33\right)$, $\left(0,27,12,34\right)$, $\left(0,28,14,35\right)$, $\left(0,29,16,36\right)$, $\left(0,30,18,37\right)$, $\left(0,31,20,38\right)$, $\left(0,32,22,39\right)$, $\left(0,8,2,9\right)$, $\left(0,3,5\right)$, $\left(0,46,42\right)$, $\left(0,44,43\right)$\\
\hline
$10$ & $\left(0,25,1,43\right)$, $\left(0,44,4,45\right)$, $\left(0,26,10,33\right)$, $\left(0,27,12,34\right)$, $\left(0,28,14,35\right)$, $\left(0,29,16,36\right)$, $\left(0,30,18,37\right)$, $\left(0,31,20,38\right)$, $\left(0,32,22,39\right)$, $\left(0,8,2,9\right)$, $\left(0,3,5\right)$, $\left(0,50,46\right)$, $\left(0,48,47\right)$\\
\hline
$11$ & $\left(0,30,13,38\right)$, $\left(0,31,15,39\right)$, $\left(0,32,17,40\right)$, $\left(0,46,18,47\right)$, $\left(0,33,19,41\right)$, $\left(0,34,21,42\right)$, $\left(0,48,22,49\right)$, $\left(0,35,23,43\right)$, $\left(0,36,25,44\right)$, $\left(0,37,27,45\right)$, $\left(0,8,2,9\right)$, $\left(0,3,5\right)$,$\left(0,54,50\right)$, $\left(0,52,51\right)$\\
\hline
$12$ & $\left(0,36,26,41\right)$, $\left(0,35,24,40\right)$, $\left(0,34,22,39\right)$, $\left(0,33,20,38\right)$, $\left(0,32,18,37\right)$, $\left(0,47,27,53\right)$, $\left(0,46,25,52\right)$, $\left(0,45,23,51\right)$, $\left(0,44,21,50\right)$, $\left(0,43,19,49\right)$, $\left(0,42,17,48\right)$, $\left(0,8,2,9\right)$, $\left(0,3,5\right)$, $\left(0,58,54\right)$,$\left(0,56,55\right)$\\
\hline
$13$ &  $\left(0,40,26,45\right)$, $\left(0,39,24,44\right)$, $\left(0,38,22,43\right)$, $\left(0,37,20,42\right)$, $\left(0,36,18,41\right)$, $\left(0,51,27,57\right)$, $\left(0,50,25,56\right)$, $\left(0,49,23,55\right)$, $\left(0,48,21,54\right)$, $\left(0,47,19,53\right)$, $\left(0,46,17,52\right)$, $\left(0,11,1,13\right)$, $\left(0,8,2,9\right)$, $\left(0,3,5\right)$, $\left(0,59,60\right)$, $\left(0,58,62\right)$\\
\hline
$14$ & $\left(0,55,1,59\right)$, $\left(0,56,4,57\right)$, $\left(0,60,10,61\right)$, $\left(0,30,13,38\right)$, $\left(0,31,15,39\right)$, $\left(0,32,17,40\right)$, $\left(0,46,18,47\right)$, $\left(0,33,19,41\right)$, $\left(0,34,21,42\right)$, $\left(0,48,22,49\right)$, $\left(0,35,23,43\right)$, $\left(0,36,25,44\right)$, $\left(0,37,27,45\right)$, $\left(0,8,2,9\right)$, $\left(0,3,5\right)$, $\left(0,66,62\right)$, $\left(0,64,63\right)$\\
\hline
$15$ & $\left(0,57,1,61\right)$, $\left(0,58,4,59\right)$, $\left(0,62,10,63\right)$, $\left(0,30,13,38\right)$, $\left(0,64,14,65\right)$, $\left(0,31,15,39\right)$, $\left(0,32,17,40\right)$, $\left(0,46,18,47\right)$, $\left(0,33,19,41\right)$, $\left(0,34,21,42\right)$, $\left(0,48,22,49\right)$, $\left(0,35,23,43\right)$, $\left(0,36,25,44\right)$, $\left(0,37,27,45\right)$, $\left(0,8,2,9\right)$, $\left(0,3,5\right)$, $\left(0,70,66\right)$, $\left(0,68,67\right)$\\
\hline
$16$ & $\left(0,46,36,53\right)$, $\left(0,45,34,52\right)$, $\left(0,44,32,51\right)$, $\left(0,43,30,50\right)$, $\left(0,42,28,49\right)$, $\left(0,41,26,48\right)$, $\left(0,40,24,47\right)$, $\left(0,61,37,69\right)$, $\left(0,60,35,68\right)$, $\left(0,59,33,67\right)$, $\left(0,58,31,66\right)$, $\left(0,57,29,65\right)$, $\left(0,56,27,64\right)$, $\left(0,55,25,63\right)$, $\left(0,54,23,62\right)$, $\left(0,8,2,9\right)$, $\left(0,3,5\right)$, $\left(0,71,72\right)$, $\left(0,70,74\right)$\\
\hline
$17$ & $\left(0,50,36,57\right)$, $\left(0,49,34,56\right)$, $\left(0,48,32,55\right)$, $\left(0,47,30,54\right)$, $\left(0,46,28,53\right)$, $\left(0,45,26,52\right)$, $\left(0,44,24,51\right)$, $\left(0,65,37,73\right)$, $\left(0,64,35,72\right)$, $\left(0,63,33,71\right)$, $\left(0,62,31,70\right)$, $\left(0,61,29,69\right)$, $\left(0,60,27,68\right)$, $\left(0,59,25,67\right)$, $\left(0,58,23,66\right)$, $\left(0,11,1,13\right)$, $\left(0,8,2,9\right)$, $\left(0,3,5\right)$, $\left(0,78,74\right)$, $\left(0,76,75\right)$\\
\hline
$18$ & $\left(0,54,36,61\right)$, $\left(0,53,34,60\right)$, $\left(0,52,32,59\right)$, $\left(0,51,30,58\right)$, $\left(0,50,28,57\right)$, $\left(0,49,26,56\right)$, $\left(0,48,24,55\right)$, $\left(0,69,37,77\right)$, $\left(0,68,35,76\right)$, $\left(0,67,33,75\right)$, $\left(0,66,31,74\right)$, $\left(0,65,29,73\right)$, $\left(0,64,27,72\right)$, $\left(0,63,25,71\right)$, $\left(0,62,23,70\right)$, $\left(0,13,1,17\right)$, $\left(0,14,4,15\right)$, $\left(0,8,2,9\right)$, $\left(0,3,5\right)$, $\left(0,78,82\right)$, $\left(0,79,80\right)$\\
\hline
$19$ & $\left(0,58,36,65\right)$, $\left(0,57,34,64\right)$, $\left(0,56,32,63\right)$, $\left(0,55,30,62\right)$, $\left(0,54,28,61\right)$, $\left(0,53,26,60\right)$, $\left(0,52,24,59\right)$, $\left(0,73,37,81\right)$, $\left(0,72,35,80\right)$, $\left(0,71,33,79\right)$, $\left(0,70,31,78\right)$, $\left(0,69,29,77\right)$, $\left(0,68,27,76\right)$, $\left(0,67,25,75\right)$, $\left(0,66,23,74\right)$, $\left(0,15,1,19\right)$, $\left(0,16,4,17\right)$, $\left(0,20,10,21\right)$, $\left(0,8,2,9\right)$, $\left(0,3,5\right)$, $\left(0,82,86\right)$, $\left(0,84,83\right)$.\\
\hline

\end{longtable}


\begin{thebibliography}{99}
\bibitem{ahmad}A. Alkasasbeh, D. Dyer, and N. Shalaby, Applying Skolem Sequences to Gracefully Label New Families of Triangular Windmills, submitted, {\em Discrete Applied Mathematics}.

\bibitem{baker2}C. Baker, and J. D. A. Manzer, Skolem-labeling of generalized three-vane windmills, {\em Australas. J. Combin.} 41 (2008), 175--204.

\bibitem{baker}C. Baker, R. Nowakowski, N. Shalaby, and A. Sharary, $M$-fold and extended $M$-fold Skolem sequences, {\em Utilitas Mathematica}, 45 (1994), 153--167.
%\bibitem{burgess}A. Burgess, Construction of Skolem and Related Sequences, Memorial University Honours Thesis, 2003.

%%%%%%%%%%%%%%%%%%%%%%%%%%%%%%%%%%%%%%%%%%%%%%%
\bibitem{Bermond} J. C. Bermond, Graceful graphs, radio antennae and French windmills, {\em Graph Theory and Combinatorics (ed. R.J. Wilson)}, Pitman, London (1979), 18--37.
\bibitem{Bermond1}J. C. Bermond, A.E. Brouwer, A. Germa, Syst\'{e}mes de triplets et diff\'{e}rences associ\'{e}es. Probl\'{e}mes combinatoires et th\'{e}orie des graphes (Colloq. Internat. CNRS, Univ. Orsay, Orsay, 1976), pp. 35--38, Colloq. Internat. CNRS, 260, CNRS, Paris, 1978.


%\bibitem{crc} C.J. Colbourn, and J.H. Dinitz (editors), {\em The CRC Handbook of Combinatorial Designs}, CRC Press, 2007, 612--616.
\bibitem{dyer} D. Dyer, I. Payne, N. Shalaby, and B. Wicks, On the graceful conjecture for triangular cacti, {\em Australas. J. Combin.} {53} (2012), 151--170.

\bibitem{Gallian} J. A. Gallian, A Dynamic Survey of Graph Labelling, {\em Electron. J. Combin.}, 5 (1998), 6--34.
\bibitem{hilton} A. J. W. Hilton, Steiner and similar triple systems, {\em Mathematica Scandinavica} 24 (1969), 208--216.

\bibitem{kejie} M. Kejie, Gracefulness of $P\left(n_1,n_2,\ldots,n_m\right)$ and $D_{m,4}$, {\em Applied Math}, 4 (1989), 95--97.
%\bibitem{heffter1}L. Heffter, Uber tripelsystemes, {\em Math. Ann.} {49} (1897), 101--112.

\bibitem{koh}K. M. Koh, D. G. Rogers, P. Y. Lee, and C. W. Toh, On graceful graphs V: unions of graphs with one vertex in common, {\em Nanta Math.}, 12 (1979), 133--136.

%\bibitem{koh}K. M. Koh, D. G. Rogers, P. Y. Lee, and C. W. Toh, On graceful graphs V: unions of graphs with one vertex in common, Nanta Math., 12 (1979) 133-136.

%\bibitem{linek} V. Linek and Z. Jiang, Extended Langford sequences with small defects, {\em J. Combin. Theory Ser.} A 84 (1998), 38--54.
\bibitem{marsh} D. C. B. Marsh, Solution of Problem E1845, {\em Amer. Math. Monthly} 74 (1967), 591--592.

\bibitem{okeefe}E. S. O'Keefe, Verification of a conjecture of Th. Skolem, {\em Math. Scand.} {9} (1961), 80--82.

\bibitem{rosa1}A. Rosa, On certain valuations of the vertices of a graph, Theory of Graphs (Internat. Symposium, Rome, July 1966), Gordon and Breach, N. Y. and Dunod Paris (1967), 349--355.
%\bibitem{rosa} A. Rosa, Cyclic Steiner Triple Systems and Labellings of Triangular Cacti, {\em Scientia Series A: Mathematical Sciences}, 1 (1988), 87--95.
\bibitem{shalaby} N. Shalaby, The existence of near-Skolem and hooked near-Skolem sequences, {\em Discrete Math.} 135 (1994), 303--319.
\bibitem{simpson} J. E. Simpson, Langford Sequences; Perfect and Hooked, {\em Discrete Math.} 44 (1983), 97--104.
\bibitem{skolem} T. Skolem, On Certain Distributions of Integers in Pairs with Given Differences, {\em Math. Scand.} 5 (1957), 57--68.
%\bibitem{shalaby1}N. Shalaby, Skolem sequences: Generalizations and Applications, Ph.D Thesis, McMaster University, Canada, 1992.

\bibitem{yang11}X. Xu, Y. Yang, H. Li, and Y. Xi, The graphs $C_{11}^{t}$ are graceful for $t \equiv 0,1\ ({\rm mod}\ 4)$, {\em Ars Combin.} 88 (2008), 429--435.

\bibitem{yang13}X. Xu, Y. Yang, L. Han, and H. Li, The graphs $C_{13}^{t}$ are graceful for $t \equiv 0,3\ ({\rm mod}\ 4)$, {\em Ars Combin.} 90 (2009), 25--32.

\bibitem{yang}Y. Yang, X. Lin, and C. Yu, The graphs $C_{5}^{t}$ are graceful for $t \equiv 0,3\ ({\rm mod}\ 4)$, {\em Ars Combin.} 74 (2005), 239--244.
\bibitem{yang7}Y. Yang, X. Xu, Y. Xi, H. Li, and K. Haque, The graphs $C_{7}^{t}$ are graceful for $t \equiv 0,1\ ({\rm mod}\ 4)$, {\em Ars Combin.} 79 (2006), 295--301.
\bibitem{yang9}Y. Yang, X. Xu, Y. Xi, and H. Huijun, The graphs $C_{9}^{t}$ are graceful for $t \equiv 0,3\ ({\rm mod}\ 4)$, {\em Ars Combin.} 85 (2007), 361--368.

%\bibitem{simpson} J.E. Simpson, Langford Sequences; Perfect and Hooked, {\em Discrete Mathematics} 44 (1983), 97--104.
%\bibitem{skolem} T. Skolem, On Certain Distributions of Integers in Pairs with Given Differences, {\em Math. Scand.} 5 (1957), 57--68.
%\bibitem{okeefe}E. S. O'Keefe, Verification of a conjecture of Th. Skolem, {\em Math. Scand.} {9} (1961), 80-82.




\end{thebibliography}
\end{document}